\documentclass[ieot]{birkmult}

\usepackage{amsmath,amsfonts}
\usepackage{amsthm}
\usepackage{enumerate}

\usepackage[dvips]{color}

\newtheorem{theorem}{Theorem}[section]
\newtheorem{proposition}[theorem]{Proposition}
\newtheorem{lemma}[theorem]{Lemma}
\newtheorem{corollary}[theorem]{Corollary}

\theoremstyle{definition}
\newtheorem{definition}[theorem]{Definition}
\newtheorem{example}[theorem]{Example}

\newtheorem{condition}[theorem]{Condition}

\newtheorem{remark}[theorem]{Remark}

\newtheorem*{remark*}{Remark}

\numberwithin{equation}{section}

\newcommand{\nC}{\mathbb C}

\newcommand{\nR}{\mathbb R}

\newcommand{\cD}{{\mathcal D}}

\newcommand{\cF}{{\mathcal F}}

\newcommand{\cK}{{\mathcal K}}
\newcommand{\cL}{{\mathcal L}}

\newcommand{\kip}{[\,\cdot\, , \cdot\,]}

\newcommand{\la}{\langle}
\newcommand{\ra}{\rangle}
\newcommand{\bpm}{\begin{bmatrix}}
\newcommand{\epm}{\end{bmatrix}}

\newcommand{\w}[1]{{#1}}

\newcommand{\co}[1]{\overline{#1}}
\newcommand{\V}[1]{\mathbf{#1}}

\newcommand{\M}[1]{\mathsf{#1}}

\DeclareMathOperator{\sgn}{sgn} 
\DeclareMathOperator{\dom}{\cD} \DeclareMathOperator{\fdom}{\cF}

\begin{document}

\title[Riesz bases of root vectors, II]%
{Riesz bases of root vectors of \\
indefinite Sturm-Liouville problems \\
with  eigenparameter  dependent \\
boundary conditions. II}

\author{Paul Binding}
\address{Department of Mathematics and Statistics, \\ University of
Calgary, \\ Calgary,  Alberta, Canada, T2N 1N4}
 \email{binding@ucalgary.ca}
\author{Branko \'{C}urgus}
\address{Department of Mathematics, \\ Western Washington University,
\\ Bellingham, WA 98225, USA}
 \email{curgus@cc.wwu.edu}


% \subjclass
% 34B09, 34B25, 47B25

% \subjclass{Primary 47B50, 47B25}


\date{\today}

\begin{abstract}
We consider a regular indefinite Sturm-Liouville problem with two
self-adjoint boundary conditions affinely dependent on the
eigenparameter.  We give sufficient conditions under which the root
vectors of this Sturm-Liouville problem can be selected to form a
Riesz basis of a corresponding weighted Hilbert space.
\end{abstract}

\maketitle

\section{Introduction} \label{s2}
Consider the following eigenvalue problem
\begin{align*}
 -f''(x) & = \lambda \, (\sgn x) f(x), \quad x \in [-1,1], \\
f'(1) & = \lambda  \, f(-1), \\
-f'(-1) & = \lambda \, f(1).
\end{align*}
Lengthy but straightforward calculations show the following: there
exist an infinite number of simple nonzero eigenvalues which
accumulate at both $-\infty$ and $+\infty$; the number $0$ is also an
eigenvalue with geometric multiplicity $1$ and algebraic multiplicity
$2$. Details can be found at the second author's web-site. It is
natural to consider this problem in the Hilbert space
$L_{2}(-1,1)\oplus \nC^2$. To our knowledge the following related
question, which presents interesting mathematical challenges, has not
been addressed.  Is it possible to select eigenvectors of the
eigenvalue problem to form a Riesz basis of the above Hilbert space?
In this article we answer such questions for a wide class of
indefinite Sturm-Liouville problems with $\lambda$-dependent boundary
conditions. In particular, our Theorem~\ref{tn=k2} applies to the
above simple example.

We consider a regular indefinite Sturm-Liouville eigenvalue problem
of the form
\begin{equation} \label{sl1}
-(p\,f')'+q\,f = \lambda \, r\, f \ \ \ \ \text{on} \ \ \ \
[-1,1].
\end{equation}
We assume throughout that the coefficients $1/p, q, r$ in
\eqref{sl1} are real and integrable over $[-1,1]$, $p(x) > 0$, and
$x\, r(x)
> 0$ for almost all $x \in [-1,1]$.  We impose the following
eigenparameter dependent boundary conditions on equation \eqref{sl1}:
\begin{equation} \label{bc1}
\M{M} \V{b}(f) = \lambda\, \M{N} \V{b}(f) ,
\end{equation}
where $\M{M}$ and $\M{N}$ are $2 \times 4$ matrices and the boundary
mapping $\V{b}$ is defined for all $f$ in the domain of \eqref{sl1}
by
\begin{equation*}
\V{b}(f) = \begin{bmatrix} f(-1) & f(1) & (pf')(-1) &
(pf')(1)\end{bmatrix}^T.
\end{equation*}
For our opening example
\[
\M{M}  =  \begin{bmatrix}
0 & 0 & 0 & 1 \\[5pt]
0 & 0 & -1 & 0
\end{bmatrix}, \quad \quad   \M{N} = \begin{bmatrix}
1 & 0 & 0 & 0 \\[5pt]
0 & 1 & 0& 0
\end{bmatrix}.
\]

We remark that more general boundary conditions have been studied by
many authors, recently for example in \cite{BJ} and \cite{BT}, but
expansion theorems were not considered. Expansion theorems for
polynomial boundary conditions and more general operators, but with
weight $r = 1$, were given in \cite{D} and \cite{Ru}.

In this article we study the problem \eqref{sl1}, \eqref{bc1} in an
operator theoretic setting established in \cite{BC}. Under
Condition~\ref{cbc} below, a definitizable self-adjoint operator $A$
in the Krein space $L_{2,r}(-1,1)\oplus \nC^2_{\Delta}$ (actually $A$
is quasi-uniformly positive as defined in \cite{CN}) is associated
with the eigenvalue problem \eqref{sl1}, \eqref{bc1}. Here $\Delta$
is a $2\times 2$ nonsingular Hermitean matrix which is determined by
$\M{M}$ and $\M{N}$; see Section~\ref{sop} for details. We remark
that the topology of this Krein space is that of the corresponding
Hilbert space $L_{2,|r|}\oplus \nC^2_{|\Delta|}$. Here, and in the
rest of the paper, we abbreviate $L_{2,r}(-1,1)$ to $L_{2,r}$ and
$L_{2,|r|}(-1,1)$ to $L_{2,|r|}$. For more details about Krein spaces
and their operators see the standard reference \cite{LaD} and
\cite{ABT} for recent developments.

Our main goal in this paper is to provide sufficient conditions on
the coefficients in \eqref{sl1}, \eqref{bc1} under which there is a
Riesz basis of the above Hilbert space consisting of the union of
bases for all the root subspaces of the above operator $A$. This will
be referred to for the remainder of this section as the {\em Riesz
basis property of} $A$. We remark that the Riesz basis property of
$A$ is equivalent, modulo a finite dimensional subspace, to
similarity of $A$ to a self-adjoint operator in a Hilbert space. The
latter similarity has been the subject of several recent papers (see
for example \cite{KKM} and \cite{KM}) involving Sturm-Liouville
expressions on $\nR$ without boundary conditions.

Existence of Riesz bases and expansion theorems with a stronger
topology, but in a smaller space corresponding to the form domain of
the operator $A$ (which in our case is a Pontryagin space), have been
considered by many authors; see \cite{BC, T} and the references
there. The results in \cite{BC} turned out to be independent of the
number and the nature of the boundary conditions and the coefficients
$p$ and $r$. In contrast, the Riesz basis property depends
nontrivially on the problem data even for the case when the boundary
conditions are $\lambda$-independent (corresponding to $\M{N} =
\M{0}$ in our notation).

Sufficient conditions on $r$ (near the turning point $0$) for the
Riesz basis property when $\M{N} = \M{0}$ can be found in \cite{Be,
CL, F3, FVY, Py, Py2}, for example. That some condition is necessary,
even in the case $p = 1$, was shown by Volkmer \cite{V} who proved
the existence of an odd $r$ for which the Dirichlet problem
\eqref{sl1} does not have this property. Recently Parfenov \cite{P}
gave a necessary and sufficient condition on an odd weight function
$r$, near its turning point $0$, for the Dirichlet problem
\eqref{sl1} to have the Riesz basis property. In \cite{BC1} we
constructed an odd $r$ for which the Dirichlet problem \eqref{sl1}
has the Riesz basis property but the anti-periodic problem does not.
This example shows that an additional condition on $r$ near the
boundary of $[-1,1]$ (which in some cases behaves as a second turning
point, in addition to $0$, for \eqref{sl1}) is needed for the general
case of \eqref{bc1}. Such conditions are given in \cite{CL} for
$\lambda$-independent boundary conditions and in \cite{BC2} for
exactly one $\lambda$-dependent boundary condition (i.e., when
$\M{N}$ has rank $1$).

In this paper we consider the more difficult case of two
$\lambda$-dependent boundary conditions. The method we use has its
origins in the work of Beals \cite{Be}. Subsequently it was developed
in \cite{C} into a criterion (given below as Theorem~\ref{W})
equivalent to the Riesz basis property of $A$.  This criterion
involves a positive homeomorphism $W$ of the Krein space
$L_{2,r}\oplus\nC^2_{\Delta}$ with the form domain of $A$ as an
invariant subspace. The explicit description of the form domain of
$A$ (given in Section \ref{sop}) depends entirely on the number $k
\in \{0,1,2\}$ of boundary conditions which do not include
derivatives in the $\lambda$-terms. We call such boundary conditions
{\em essential}. Note that this differs from the usual terminology
for $\lambda$-independent conditions.  For example, in our
terminology $y'(1) = \lambda y(1)$ is an essential boundary
condition.

The direct sum structure of the Krein space
$L_{2,r}\oplus\nC^2_{\Delta}$ naturally leads us to consider the
homeomorphism $W$ as a block operator matrix, the top left entry
$W_{11}$ being an operator on $L_{2,r}$. Since it is clear from
Section \ref{sop} that the functional components of the vectors in
the form domain of $A$ are (absolutely) continuous, we see that
$W_{11}$ induces a boundary matrix $\M{B}$ satisfying
\[
\M{B} \begin{bmatrix}  f(-1) \\ f(1)
\end{bmatrix} = \begin{bmatrix}  (W_{11}f)(-1) \\ (W_{11}f)(1)
\end{bmatrix}.
\]
An important hurdle, with analogues in several of above references,
is to solve the inverse problem of finding a suitable $W_{11}$ for a
given matrix $B$. For example, in \cite{BC2} (see also
Section~\ref{s0-11} below) such operators $W_{11}$ were constructed
with special diagonal $\M{B}$ under one-sided Beals type conditions
at $-1$ or $1$.   In Section~\ref{smpm1} we use conditions at $-1$,
at $1$, and a condition connecting $-1$ and $1$ to produce $W_{11}$
with an arbitrary prescribed boundary matrix $\M{B}$.

In Sections~\ref{sk=0} and \ref{sn2k1} we complete the construction
of $W$, thus establishing our sufficient conditions for the Riesz
basis property.  When there are no essential boundary conditions
($k=0$), it turns out that the one-sided Beals type condition at $0$
suffices; see Theorem~\ref{tk0sebc}. In other cases, however, we
need conditions near the boundary of $[-1,1]$.  Conditions at $0$,
and at $-1$ or $1$, are sufficient if $k=2$ and $\Delta$ is
definite. If $\Delta$ is indefinite, then we also need the condition
linking $-1$ and $1$. In these cases it suffices to construct $W$ as
a block diagonal matrix. This is carried out in Theorem~\ref{tn=k2}.

The most difficult case is $k=1$ which we tackle in
Section~\ref{sn2k1}.  In this case we need not only off-diagonal
blocks for $W$, but also a perturbation $K$ of $W_{11}$, where $K$
is an integral operator whose construction is rather delicate.  Our
final result Theorem~\ref{tn2k1} is as follows. If only one boundary
point $-1$ or $1$ appears with $\lambda$ in the essential boundary
condition, then a Beals type condition at that point and at $0$ are
sufficient. Otherwise we need conditions at both boundary points and
at $0$, as well as the condition linking $-1$ and $1$.

To conclude this introduction we remark that our conditions simplify
drastically if $p$ is even and $r$ is odd, a case which has been
studied by several authors \cite{BC1, P, V}. In fact all the
conditions that we impose on the boundary are then equivalent; see
Example~\ref{r42} and Corollary~\ref{cpero}.

\section{Operators associated with the eigenvalue problem}
\label{sop}

The maximal operator $S_{\max}$ in $L_{2,r}$ associated with
\eqref{sl1} is defined by
\[
S_{\max} : f \mapsto \ell(f) := \frac{1}{r} \bigl(-(pf')' + qf
\bigr), \ \ \ f \in \cD(S_{\max}),
\]
where
\[
\dom(S_{\max}) =  \cD_{\max} = \bigl\{ f \in L_{2,r} : f, pf' \in
AC[0,1], \ \ell(f) \in L_{2,r} \bigr\}.
\]
We define the boundary mapping $\V{b}$ by
\begin{equation*}
\V{b}(f) = \begin{bmatrix} f(-1) & f(1) & (pf')(-1) &
(pf')(1)\end{bmatrix}^T, \ \ \  f \in \dom(S_{\max}).
\end{equation*}
and the concomitant matrix $\M{Q}$ corresponding to $\V{b}$ by
\[
\M{Q} = \text{{\Large $i$}} \begin{bmatrix}
  0  &  0 &  - 1 & 0  \\[2pt]
  0 & 0 & 0 & 1  \\[2pt]
  1 & 0 & 0 & 0  \\[2pt]
  0 & -1  & 0 & 0
\end{bmatrix}.
\]
The significance of $\M{Q}$ is captured by the following identity
\[
\int_{-1}^1 \bigl( S_{\max} f \overline{g} - f S_{\max} \overline{g}
\bigr) r = i \, \V{b}(g)^* \M{Q} \V{b}(f), \ \ f,g \in \cD_{\max}.
\]
We note that $\M{Q} = \M{Q}^{-1}$.

Throughout, we shall impose the following nondegeneracy and
self-adjointness condition on the boundary data.

\begin{condition} \label{cbc}
The boundary matrices $\M{M}$ and $\M{N}$ in \eqref{bc1} satisfy the
following:
\begin{enumerate}[(1)]
\item
the $4 \times 4$ matrix $\begin{bmatrix} \M{M} \\
\M{N}\end{bmatrix}$ is nonsingular,
\item
 $\M{M Q M}^* = \M{N Q N}^* = \M{0}$,
\item \label{ibp4}
the $2 \times 2$ matrix $i \M{M} \M{Q}^{-1} \M{N}^*$ is self-adjoint
and invertible and we define
\[
\M{\Delta} := - i \bigl(\M{M} \M{Q}^{-1} \M{N}^*\bigr)^{-1}.
\]
\end{enumerate}
\end{condition}

Clearly the boundary value problem \eqref{sl1},\eqref{bc1} will not
change if row reduction is applied to the coefficient matrix
\begin{equation} \label{cm}
\begin{bmatrix} \M{M} & \M{N}  \end{bmatrix}.
\end{equation}
In what follows we will assume that the matrix in \eqref{cm} is row
reduced to row echelon form (starting the reduction at the bottom
right corner). In particular the matrix $\M{N}$ has the form
\begin{equation*} %\label{Le}
\M{N} =  \begin{bmatrix} \M{N}_{e}  & \M{0}   \\[5pt]
   \M{N}_{1} & \M{N}_{n}    \end{bmatrix}.
\end{equation*}
The matrix $\M{0}$ in the formula for $\M{N}$ is $k \times 2$ with
$k \in \{0,1,2\}$.  The $k \times 2$ matrix $\M{N}_{e}$ and the
$(2-k)\times 2$ matrix $\M{N}_{n}$ are of maximal ranks.

There are three possible cases for $\M{N}$ in \eqref{cm}:
\begin{enumerate}[(a)]
\item
\ $\M{N}_n$ is a $2 \times 2$ identity matrix (so $k = 0$),
\item
\ $\M{N}_e$ and $\M{N}_n$ are nonsingular $1 \times 2$ (row)
matrices  (so $k = 1$),
\item
\ $\M{N}_e$ is a $2 \times 2$ identity matrix  (so $k = 2$).
\end{enumerate}
In case (a), both boundary conditions in \eqref{bc1} are {\em
non-essential}, that is both rows on the right hand side of
\eqref{bc1} contain derivatives. In case (b), the boundary condition
corresponding to the first row in \eqref{bc1} is {\em essential},
that is no derivatives appear in this row on the right hand side;
the second boundary condition in \eqref{bc1} is non-essential. In
case (c), both boundary conditions in \eqref{bc1} are essential.
Evidently $k$ is the number of essential boundary conditions.

Next we define a Krein space operator associated with the problem
\eqref{sl1},\eqref{bc1}.  We consider the linear space $L_{2,r}
\oplus \nC_{\Delta}^2$, equipped with the inner product
\[
\left[ \begin{pmatrix} f \\ \V{u} \end{pmatrix} ,
\begin{pmatrix} g
\\ \V{v} \end{pmatrix} \right] := \int_{-1}^{1} f \co{g} r \,
+ \, \V{v}^{*} \Delta \V{u}, \ \ \  f, g \in L_{2,r}, \
\V{u},\V{v} \in \nC^2.
\]
Then $\bigl(L_{2,r} \oplus \nC_{\Delta}^2,\kip\bigr)$ is a Krein
space.  A fundamental symmetry on this Krein space is given by
\begin{equation*} %\label{fsy}
J := \begin{bmatrix} J_0  & 0  \\[5pt]  0 & \sgn(\M{\Delta})
\end{bmatrix},
\end{equation*}
where $2 \times 2$ matrix $\sgn(\M{\Delta})$ and $J_0:L_{2,r} \to
L_{2,r}$ are defined by
\[
\sgn(\Delta) = |\Delta|^{-1} \Delta \ \ \ \text{and} \ \ \
 (J_0f)(t) :=  f(t)\, \sgn(r(t)), \ \ t \in [-1,1].
\]
Then $\la\,\cdot\,,\,\cdot\,\ra := [J\,\cdot\,,\,\cdot\,]$ is a
positive definite inner product which turns $L_{2,r} \oplus
\nC_{\Delta}^2$ into a Hilbert space $\bigl( L_{2,|r|} \oplus
\nC_{|\Delta|}^2,\la\,\cdot\,,\,\cdot\,\ra \bigr)$.  The topology of
$L_{2,r} \oplus \nC_{\Delta}^2$ is defined to be that of $L_{2,|r|}
\oplus \nC_{|\Delta|}^2$, and a {\em Riesz basis} of $L_{2,r} \oplus
\nC_{\Delta}^2$ is defined as a homeomorphic image of an orthonormal
basis of $L_{2,|r|} \oplus \nC^2_{|\Delta|}$.

We define the operator $A$ in the Krein space $L_{2,r} \oplus
\nC_{\Delta}^2$ on the domain
\begin{equation*} % \label{domwA}
\dom(\w{A}) =  \left\{ \begin{bmatrix} f \\[3pt] \M{N} \V{b}(f)
\end{bmatrix} \in \w{\cK} \ : \  f \in \dom\bigl(S_{\max}\bigr)  \right\}
\end{equation*}
by
\begin{equation*} % \label{defwA}
\w{A}  \begin{bmatrix} f \\[3pt] \M{N} \V{b}(f)
\end{bmatrix} :=
\begin{bmatrix} S_{\max}f \\[3pt] \M{M} \V{b}(f) \end{bmatrix}, \ \
\ \ \ \ \ \ f \in \dom(\w{A}).
 \end{equation*}

Using \cite[Theorems 3.3 and 4.1]{BC} we see that this operator is
definitizable with discrete spectrum in the Krein space $L_{2,r}
\oplus \nC_{\Delta}^2$.  As in \cite[Theorem 2.2]{BC2}, we then
obtain the following, which is our basic tool.

\begin{theorem} \label{W}
Let $\fdom(A)$ denote the form domain of $A$. Then there exists a
Riesz basis of $L_{2,r} \oplus \nC^2_{\Delta}$ which consists of
root vectors of $A$ if and only if there exists a bounded,
boundedly invertible, positive operator $W$ in $L_{2,r} \oplus
\nC^2_{\Delta}$ such that
\[
W \fdom(A) \subset \fdom(A).
\]
\end{theorem}

In order to apply this result, we need to characterize the form
domain $\fdom(A)$.  To this end, let $\cF_{\max}$ be the set of all
functions $f$ in $L_{2,r}$ which are absolutely continuous on
$[-1,1]$ and such that $\int_{-1}^1 p\, |f^{\prime}|^2 < +\infty$.

By \cite[Theorem 4.2]{BC}, there are three possible cases for the
form domain $\fdom(A)$ of $A$, corresponding to cases (a), (b) and
(c) above.
\begin{enumerate}[(a)]
\item \
If $\M{N}_n =
\begin{bmatrix} 1 & 0
\\[2pt] 0 & 1 \end{bmatrix}$, then
\begin{align}
\fdom(A) & = \left\{ \begin{bmatrix} f \\[3pt] \V{v} \end{bmatrix} \in
 \begin{matrix} L_{2,r} \\ \oplus \\ \nC^2_{\Delta} \end{matrix}
  \ : \ f \in \cF_{\max}, \  \V{v} \in \nC^2
\right\}. \label{fd1}
 \intertext{
\item \
If $\M{N}_e  = [u \ \ v]$ with $u,v \in \nC$ and $|u|^2 + |v|^2 \neq
0$, then }
 \fdom(\w{A}) & = \left\{ \begin{bmatrix}
  f \\[2pt]
  u f(-1) + v f(1) \\[2pt]
  z  \end{bmatrix} \in
  \begin{matrix} L_{2,r} \\ \oplus \\ \nC^2_{\Delta} \end{matrix}
    \ : \ f \in \cF_{\max}, \ z \in \nC \right\}.
    \nonumber % \label{fd2}
 \intertext{
\item \
If $\M{N}_e = \begin{bmatrix} 1 & 0
\\[2pt] 0 & 1 \end{bmatrix}$, then  }
 \fdom(A) & = \left\{ \begin{bmatrix} f \\[2pt] f(-1) \\[2pt] f(1)
\end{bmatrix} \in
\begin{matrix} L_{2,r} \\ \oplus \\ \nC_{\Delta}^2 \end{matrix}
 \ : \ f \in \cF_{\max} \right\}.\label{fd3}
 \end{align}
\end{enumerate}

To construct an operator $W$ as in Theorem~\ref{W} we need to impose
conditions (to be given in the next two sections) on the
coefficients $p$ and $r$ in \eqref{sl1}. In all cases we need
Condition~\ref{c0} in a neighborhood of $0$, and in some cases we
need one of two Conditions, \ref{cat-1} or \ref{cat1}, on $r$ in
neighborhoods of $-1$ or $1$. These will be discussed in
Section~\ref{s0-11}. In some cases we also need Condition~\ref{mc1}
connecting the boundary points $-1$ and $1$. This is developed in
Section~\ref{smpm1}.

\section{Conditions at $0$, $-1$ and $1$}
\label{s0-11}

In this section we recall the remaining concepts and results from
\cite[Sections~3, 4 and 5]{BC2} which we need in this paper.

A closed interval of non-zero length is said to be a {\em left
half-neighborhood} of its right endpoint and a {\em right
half-neighborhood} of its left endpoint. Let $\imath$ be a closed
subinterval of $[-1,1]$.  By $\cF_{\max}(\imath)$  we denote the
set of all functions $f$ in $L_{2,r}(\imath)$ which are absolutely
continuous on $\imath$ and such that $\int_{\imath} p\,
|f^{\prime}|^2 < +\infty$. With this notation we have $\cF_{\max}
= \cF_{\max}[-1,1]$.

\begin{definition} \label{dscab}
Let $p$ and $r$ be the coefficients in \eqref{sl1}. Let $a, b \in
[-1,1]$ and let $h_a$ and $h_b$, respectively,  be
half-neighborhoods of $a$ and $b$ which are contained in $[-1,1]$.
We say that the ordered pair $(h_a, h_b)$ is {\em smoothly
connected} if there exist
\begin{enumerate}[(a)\ \ ]
 \item
positive real numbers $\epsilon$ and $\tau$,
 \item
non-constant affine functions $\alpha:[0,\epsilon] \rightarrow h_a$
and $\beta:[0,\epsilon] \rightarrow h_b$,
 \item
non-negative real functions $\rho$ and $\varpi$ defined on
$[0,\epsilon]$
\end{enumerate}
such that
\begin{enumerate}[(i)\ \ ]
\item \label{idscabi}
$\alpha(0) = a$ \ and \ $\beta(0) = b$,
 \item \label{idscabii}
$p\circ\alpha$ and $p\circ\beta$ are locally integrable on the
interval $(0,\epsilon]$,
 \item \label{idscabiii}
$\rho\circ \alpha^{-1} \in
\cF_{\max}\bigl(\alpha([0,\epsilon])\bigr)$,
\item \label{idscabiiia}
$1/\tau < \varpi < \tau$ a.e. on $[0,\epsilon]$,
 \item \label{idscabiv}
$\displaystyle \rho(t)
  = \frac{\,\bigl|r\bigl(\beta(t)\bigr)\bigr|}%
  {\bigl|r\bigl(\alpha(t)\bigr)\bigr|}$
  \ \  and \ \
  $\displaystyle   \varpi(t)  =
\frac{\,p\bigl(\beta(t)\bigr)}{p\bigl(\alpha(t)\bigr)}$ \ \ \
 for  \ \ \ $\displaystyle  t \in (0,\epsilon]$.
\end{enumerate}
The numbers $\alpha', \beta'$ (the slopes of $\alpha$, $\beta$,
respectively) and $\rho(0)$ are called the {\em parameters} of the
smooth connection.
\end{definition}

A broad class of examples satisfying this definition can be given
via the following one.
\begin{definition}
Let $\nu$ and $a$ be real numbers and let $h_a$ be a
half-neighborhood of $a$.  Let $g$ be a function defined on $h_a$.
Then $g$ is {\em of order} $\nu$ {\em on} $h_a$ if there exists $g_1
\in C^1(h_a)$ such that
\begin{equation*}
 g(x) = |x-a|^{\nu}g_1(x) \ \ \ \ \text{and} \ \ \ \ g_1(x) \neq
 0, \ \ \ x \in h_a.
\end{equation*}
(The absolute value is missing in the corresponding definition in
\cite{BC2}).
\end{definition}

\begin{example} \label{e35}
Let $a, b \in \{-1,0,1\}$. Let $h_a$ and $h_b$ be half-neighborhoods
of $a$ and $b$, respectively, and contained in $[-1,1]$. For
simplicity assume that $p=1$. If $r$ in \eqref{sl1} has order $\nu$
($> -1$ to ensure integrability) on both half-neighborhoods $h_a$
and $h_b$ then as noted in \cite{BC2} the half-neighborhoods $h_a$
and $h_b$ are smoothly connected.  Moreover the parameters of the
smooth connection are nonzero numbers.  We remark that that $p$ can
be much more general -- see \cite[Example~3.4]{BC2}.
\end{example}

\begin{theorem} \label{tgenS}
Let $\imath$ and $\jmath$ be closed intervals, $\imath, \jmath \in
\bigl\{[-1,0],[0,1] \bigr\}$. Let $a$ be an endpoint of $\,\imath$
and let $b$ be an endpoint of $\,\jmath$.  Denote by $a_1$ and
$b_1$, respectively, the remaining endpoints.  Assume that the
half-neighborhoods $\imath$ of $a$ and $\jmath$ of $b$ are smoothly
connected with parameters $\alpha', \beta'$ and $\rho(0)$.  Then
there exists an operator
 $$
S: L_{2,|r|}(\imath)\rightarrow L_{2,|r|}(\jmath)
 $$
such that the following hold:
\begin{enumerate}[{\rm ($S$-1)\ }]
\item \label{itS1}
$S \in \cL\bigl(L_{2,|r|}(\imath),L_{2,|r|}(\jmath)\bigr), \ S^*\in
\cL\bigl(L_{2,|r|}(\jmath),L_{2,|r|}(\imath)\bigr)$;
\item \label{itS2}
$(Sf)(x) = 0, \ |x - b_1| \leq \frac{1}{2}$ for all $f \in
L_{2,|r|}(\imath)$ \ and \\
$(S^*g)(x) = 0, \ |x - a_1| \leq \frac{1}{2}$ for all $g \in
L_{2,|r|}(\jmath)$;
\item \label{itS4}
$S\cF_{\max}(\imath) \subset \cF_{\max}(\jmath)$, \ \
$S^*\cF_{\max}(\jmath) \subset \cF_{\max}(\imath)$;
\item \label{itS3}
For all $f \in \cF_{\max}(\imath)$ and all $g \in
\cF_{\max}(\jmath)$ we have
\[
\lim\limits_{\substack{y\rightarrow b \\ y \in \jmath\,}} \
(Sf)(y) \, = \, |\alpha'| \,
 \lim\limits_{\substack{x\rightarrow a\\ x \in
\imath}} \, f(x), \ \ \  \lim\limits_{\substack{x\rightarrow a\\
x \in \imath}} \ (S^*g)(x) \ = \ |\beta'|\rho(0) \,
\lim\limits_{\substack{y\rightarrow b \\ y \in \jmath\,}} \,
g(y).
\]
\end{enumerate}
\end{theorem}

\noindent This is \cite[Theorem 3.6]{BC2}.

\begin{condition}[Condition at $0$] \label{c0}
Let $p$ and $r$ be coefficients in \eqref{sl1}.  Denote by
$h_{0{-}}$ a generic left and by $h_{0{+}}$ a generic right
half-neighborhood of $0$.  We assume that at least one of the four
ordered pairs of half-neighborhoods
\begin{equation*} %\label{eq4p}
(h_{0{-}},h_{0{-}}), \ \ \ \ (h_{0{-}},h_{0{+}}), \ \ \ \
(h_{0{+}},h_{0{-}}), \ \ \ \ (h_{0{+}},h_{0{+}}),
\end{equation*}
is smoothly connected with the connection parameters $\alpha_0',
\beta_0'$ and $\rho_0(0)$ such that $|\alpha_0'| \neq |\beta_0'|
\rho_0(0)$.
\end{condition}

We note from Example~\ref{e35} that this condition is automatically
satisfied if $p=1$ and $r$ is of order $\nu$ on some
half-neighborhood of $0$.

\begin{theorem} \label{twat0}
Assume that the coefficients $p$ and $r$ satisfy
Condition~{\rm~\ref{c0}}.  Then there exists an operator
 $$
W_0:L_{2,r} \rightarrow L_{2,r}
 $$
such that the following hold:
\begin{enumerate}[{\rm (a)}]
\item
\ $W_0$ is bounded on $L_{2,|r|}$;
\item \label{itwat0b}
$J_0 W_0 > I$, in particular $W_0^{-1}$ is bounded and $W_0$ is
positive on the Krein space $L_{2,r}$;
\item \label{itwat0c}
$(W_0f)(x) = (J_0f)(x), \ \ \ \frac{1}{2} \leq |x| \leq 1, \ \ \ f
\in L_{2,r}$;
\item
\ $W_0\cF_{\max} \subset \cF_{\max}$.
\end{enumerate}
\end{theorem}

\noindent This is \cite[Theorem 4.2]{BC2}.

\begin{condition}[Condition at $-1$] \label{cat-1}
Let $p$ and $r$ be coefficients in \eqref{sl1}.  We assume that a
right half neighborhood of $-1$ is smoothly connected to a right
half neighborhood of $-1$ with the connection parameters
$\alpha_{-1}', \beta_{-1}'$ and $\rho_{-1}(0)$ such that
$|\alpha_{-1}'| \neq |\beta_{-1}'| \rho_{-1}(0)$.
\end{condition}

\begin{condition}[Condition at $1$] \label{cat1}
Let $p$ and $r$ be coefficients in \eqref{sl1}.  We assume that a
left half-neighborhood of $1$ is smoothly connected to a left
half-neighborhood of $1$ with the connection parameters
$\alpha_{+1}', \beta_{+1}'$ and $\rho_{+1}(0)$ such that
$|\alpha_{+1}'| \neq |\beta_{+1}'| \rho_{+1}(0)$.
\end{condition}

Again, we note from Example~\ref{e35} that these conditions are
automatically satisfied if $p=1$ and $r$ is of order $\nu_{-1}$ and
$\nu_{+1}$ on some half-neighborhood (in $[-1,1]$) of $-1$ and $1$,
respectively.

The following two propositions appear in \cite{BC2} as
Propositions~5.3 and~5.4, respectively.

\begin{proposition} \label{pwat-1}
Assume that the coefficients $p$ and $r$ satisfy
Condition~{\rm~\ref{cat-1}}. Let $b$ be an arbitrary complex number.
Then there exists an operator
 $$
W_{-1}:L_{2,r} \rightarrow L_{2,r}
 $$
such that the following hold:
\begin{enumerate}[{\rm (a)\ }]
\item
 $W_{-1}$ is bounded on $L_{2,|r|}$;
\item \label{ipwat-1b}
 $J_0W_{-1} > I$, in particular $(W_{-1})^{-1}$ is bounded and
$W_{-1}$ is positive on the Krein space $L_{2,r}$;
\item \label{ipwat-1c}
$(W_{-1}f)(x) = (J_0f)(x), \ \ \ - \frac{1}{2} \leq x \leq 1, \ \ \
f \in L_{2,r}$;
\item
$W_{-1}\cF_{\max} \subset \cF_{\max}[-1,0]\oplus\cF_{\max}[0,1]$;
\item
 $(W_{-1}f)(-1) = b f(-1) \ \ \ \text{for all} \ \ \ f \in
\cF_{\max}$.
\end{enumerate}
\end{proposition}

\begin{proposition} \label{pwat1}
Assume that the coefficients $p$ and $r$ satisfy
Condition~{\rm~\ref{cat1}}. Let $b$ be an arbitrary complex number.
Then there exists an operator
 $$
 W_{+1}:L_{2,r} \rightarrow L_{2,r}
 $$
such that the following hold:
\begin{enumerate}[{\rm (a)}]
\item
 $W_{+1}$ is bounded on $L_{2,|r|}$;
\item \label{ipwat1b}
 $J_0W_{+1} > I$, in particular $(W_{+1})^{-1}$ is bounded and
$W_{+1}$ is positive on the Krein space $L_{2,r}$;
\item \label{ipwat1c}
$(W_{+1}f)(x) = (J_0f)(x), \ \ \  - 1 \leq x \leq \frac{1}{2}, \ \ \
f \in L_{2,r}$;
\item
 $W_{+1}\cF_{\max} \subset
\cF_{\max}[-1,0]\oplus\cF_{\max}[0,1]$;
\item
 $(W_{+1}f)(1) = b f(1) \ \ \ \text{for all} \ \ \ f \in
\cF_{\max}$.
\end{enumerate}
\end{proposition}

\section{Mixed condition at $\pm 1$ and
associated operator} \label{smpm1}

In this section we establish analogues of the above results for a
new condition involving both endpoints of the interval $[-1,1$].

\begin{condition}[Condition at $-1,1$] \label{mc1}
Let $p$ and $r$ be the coefficients in \eqref{sl1}.  We assume that
at least one of the following three conditions is satisfied.
\begin{enumerate}[(A)]
\item \label{imc1A}
There are two smooth connections each connecting a right
half-neighborhood of $-1$ to a left half-neighborhood of $1$ with
the connection parameters $\alpha_{mj}'$, $\beta_{mj}'$ and
$\rho_{mj}(0)$, \ $j=1,2$, such that
\begin{equation} \label{eqdet1}
 \left| \begin{matrix} |\alpha_{m1}'| & |\alpha_{m2}'| \\[8pt]
 |\beta_{m1}'| \rho_{m1}(0) &  |\beta_{m2}'| \rho_{m2}(0)
\end{matrix} \right| \neq 0.
\end{equation}
\item \label{imc1B}
There are two smooth connections each connecting a left
half-neighborhood of $1$ to a right half-neigh\-bor\-hood of $-1$
with the connection parameters $\alpha_{mj}', \beta_{mj}'$ and
$\rho_{mj}(0)$, \ $j=1,2$, such that \eqref{eqdet1} holds.
\item \label{imc1C}
A right half-neigh\-bor\-hood of $-1$ is smoothly connected to a
left half-neigh\-bor\-hood of $1$ with the connection parameters
$\alpha_{m1}', \beta_{m1}'$ and $\rho_{m1}(0)$, and a left
half-neigh\-bor\-hood of $1$ is smoothly connected to a right
half-neigh\-bor\-hood of $-1$ with the connection parameters
$\alpha_{m2}', \beta_{m2}'$ and $\rho_{m2}(0)$, such that
\begin{equation*}  % \label{eqdet2}
 \left| \begin{matrix} |\alpha_{m1}'| &   |\beta_{m2}'| \rho_{m2}(0)
   \\[8pt]
 |\beta_{m1}'| \rho_{m1}(0) & |\alpha_{m2}'|
\end{matrix} \right| \neq 0.
\end{equation*}
\end{enumerate}
\end{condition}

\begin{example} \label{e42}
From Example~\ref{e35} it follows that this condition is satisfied
if $p=1$ and $r$ has the same order $\nu$ on a right
half-neighborhood of $-1$ and a left half-neighborhood of $1$.
\end{example}

\begin{example} \label{r42}
If $p$ is an even function and $r$ is odd, then it turns out that
Conditions~\ref{cat-1}, \ref{cat1} and \ref{mc1} are equivalent. The
first equivalence is clear. For the second, assume that
Condition~\ref{cat1} is satisfied. Let $\alpha_{+1}$ and
$\beta_{+1}$ be the corresponding affine functions from
Definition~\ref{dscab} defined on $[0,\epsilon]$. Now define
$\alpha_{m1}(t) = \alpha_{+1}(t), \beta_{m1}(t) = -\beta_{+1}(t), \,
t\in [0,\epsilon)$, so $\rho_{m1} = \rho_{+1}$. Note that $p$ is
locally integrable on $[\alpha_{+1}(\epsilon),1)$ by
Definition~\ref{dscab} (ii).  Then define $\alpha_{m2}(t) = 1 - t,
\beta_{m2}(t) = -1+t, \, t\in \bigl[0,1-\alpha_{+1}(\epsilon)\bigr)$
and so $\rho_{m2} = 1$. Then Condition~\ref{mc1}(B) is satisfied
since \eqref{eqdet1} takes the form
\begin{equation*}
 \left| \begin{matrix} |\alpha_{m1}'| & |\alpha_{m2}'| \\[8pt]
 |\beta_{m1}'| \rho_{m1}(0) &  |\beta_{m2}'| \rho_{m2}(0)
\end{matrix} \right|
 =
 \left| \begin{matrix} |\alpha_{+1}'| & 1 \\[8pt]
 |\beta_{+1}'| \rho_{+1}(0) &  1
\end{matrix} \right|
\end{equation*}
which is nonzero by Condition~\ref{cat1}. The proof of the converse
is similar.
\end{example}

\begin{example} \label{e43}
We call a function $g:[-1,1] \rightarrow \nC$ {\em nearly odd} ({\em
nearly even}) if there exists a positive constant $c \neq 1$ such
that $g(-x) = - c\, g(x)$ ($g(-x) = c\, g(x)$) for almost all $x \in
(0,1]$. We note that if $p$ is a nearly even function and $r$ is
nearly odd, both Conditions~\ref{c0} and \ref{mc1} are satisfied.
Also, Conditions~\ref{cat-1} and~\ref{cat1} are equivalent. The
verification is straightforward.
\end{example}

\begin{example} \label{e44}
Let $p=1$ and $r(x) = -1$ for $x \in [-1,0)$ and $r(x) = 1-x$ for $x
\in [0,1]$. It is not difficult to verify directly that these
functions satisfy Conditions~\ref{c0}, \ref{cat-1} and~\ref{cat1},
but not Condition~\ref{mc1}. In addition notice that $r$ is of order
$0$ in a right half-neighborhood of $-1$ and of order $1$ in a left
half-neighborhood of $1$.
\end{example}

The proof of the following theorem occupies the remainder of this
section.

\begin{theorem} \label{twac1}
Assume that the coefficients $p$ and $r$ satisfy
Conditions~{\rm~\ref{cat-1}}, {\rm \ref{cat1}} and~{\rm \ref{mc1}}.
Let $b_{jk}, \, j,k=1,2$, be arbitrary complex numbers. Then there
exists an operator
 $$
W_{s1}:L_{2,r} \rightarrow L_{2,r}
 $$
such that the following hold:
\begin{enumerate}[{\rm (a)}]
 \item \label{itwac1a}
$W_{s1}$ is bounded on the Hilbert space $L_{2,|r|}$;
 \item \label{itwac1b}
$J_0W_{s1} > I$,  in particular $W_{s1}^{-1}$ is bounded and
$W_{s1}$ is positive on the Krein space $L_{2,r}$;
 \item \label{itwac1c}
$(W_{s1}f)(x) = (J_0f)(x), \ \ \ -\frac{1}{2} \leq x \leq
\frac{1}{2}, \ \ \ f\in L_{2,r}$;
 \item \label{itwac1d}
$W_{s1} \cF_{\max} \subset \cF_{\max}[-1,0]\oplus \cF_{\max}[0,1]$;
\item \label{itwac1e} \ \vspace*{-7pt}
\[
\begin{bmatrix}
(W_{s1}f)(-1) \\[5pt]  \rule{0pt}{0pt} (W_{s1}f)(1)
\end{bmatrix} =  {\begin{bmatrix}
{b}_{11} & {b}_{12} \\[5pt]
{b}_{21} & {b}_{22}
\end{bmatrix}} \begin{bmatrix}
f(-1) \\[5pt]  \rule{0pt}{0pt} f(1)
\end{bmatrix}.
\]
\end{enumerate}
\end{theorem}
\begin{proof}
We construct $W_{s1}$ in the form
\[
W_{s1} = J_0 \bigl( X_{s1}^*X_{s1} + I \bigr),
\]
where
\[
X_{s1} = \begin{bmatrix} X_{11} & X_{12} \\[3pt] X_{21} & X_{22}
\end{bmatrix}
\]
is a block operator matrix corresponding to the decomposition
\[
L_{2,|r|}=L_{2,|r|}(-1,0)\oplus L_{2,|r|}(0,1).
\]

We split the proof into three parts. The off-diagonal and diagonal
entries of $X_{s1}$ are constructed in the first and second parts,
respectively. In the third part we establish the stated properties
of $W_{s1}$.

\medskip

\noindent{\bf 1.} \ To construct the off-diagonal operators we treat
each case (A), (B), (C) of Condition~\ref{mc1} separately.

\medskip

%
% case 1
%
\noindent {\bf Case} (A). \ By Theorem~\ref{tgenS} there exist
operators
 $$
S_{mj} : L_{2,|r|}(-1,0) \rightarrow L_{2,|r|}(0,1), \ j=1,2,
 $$
which satisfy ($S$-\ref{itS1})-($S$-\ref{itS3}) in
Theorem~\ref{tgenS} with $\imath = [-1,0],\, \jmath = [0,1]$, $a =
-1$ and $b = 1$. In particular, for $f \in \cF_{\max}[-1,0]$ and
$j=1,2$,
\begin{equation*} % \label{eqpS03m1}
(S_{mj}f)(1) = |\alpha_{mj}'| \, f(-1), \ \ \ \ \ (S_{mj}^*f)(-1) =
|\beta_{mj}'|\, \rho_{mj}(0) \, f(1).
\end{equation*}
To simplify the formulas we use the following notation
 \[
\Upsilon : =
 \left| \begin{matrix} |\alpha_{m1}'| & |\alpha_{m2}'| \\[8pt]
 |\beta_{m1}'|\, \rho_{m1}(0) & |\beta_{m2}'|\, \rho_{m2}(0)
 \end{matrix} \right|.
 \]
Define
 $$
X_{21} : L_{2,|r|}(-1,0) \rightarrow L_{2,|r|}(0,1),
 $$
by
\begin{align*}
X_{21} &
 :=  b_{21} \,  \Upsilon^{-1} \
 \left| \begin{matrix} S_{m1} & S_{m2} \\[8pt]
|\beta_{m1}'|\, \rho_{m1}(0) & |\beta_{m2}'|\, \rho_{m2}(0)
\end{matrix} \right| .
\end{align*}
Here and below we write such determinants as abbreviations
 for corresponding linear combinations of operators.
For all $f \in \cF_{\max}[-1,0]$ we have
\begin{align*}
(X_{21}f)(1) & =  b_{21} \,
 \Upsilon^{-1} \
  \left| \begin{matrix} |\alpha_{m1}'|\,f(-1) & |\alpha_{m2}'|\,f(-1)
   \\[8pt]
 |\beta_{m1}'|\, \rho_{m1}(0) & |\beta_{m2}'|\, \rho_{m2}(0)
 \end{matrix} \right|  = b_{21} f(-1).
\end{align*}
Also for all $g \in \cF_{\max}[0,1]$ we have
\begin{align*}
(X_{21}^*g)(-1) & = \co{b}_{21} \,
  \Upsilon^{-1} \
  \left| \begin{matrix}
  |\beta_{m1}'|\, \rho_{m1}(0)\,g(1)
   & |\beta_{m2}'|\, \rho_{m2}(0)\,g(1) \\[8pt]
  |\beta_{m1}'|\, \rho_{m1}(0) & |\beta_{m2}'|\, \rho_{m2}(0)
 \end{matrix} \right| = 0 .
\end{align*}

Now define the opposite off diagonal corner
$$
X_{12} : L_{2,|r|}(0,1) \rightarrow L_{2,|r|}(-1,0),
 $$
by
\begin{align*}
X_{12} & :=  -b_{12} \,
 \Upsilon^{-1} \
 \bigl( - |\alpha_{m2}'|\ S_{m1}^* +
  |\alpha_{m1}'|\ S_{m2}^* \bigr)
  = -b_{12} \,
  \Upsilon^{-1} \
  \left| \begin{matrix}
  |\alpha_{m1}'| & |\alpha_{m2}'| \\[8pt]
  S_{m1}^* & S_{m2}^*
 \end{matrix} \right| .
\end{align*}
Then for all $f \in \cF_{\max}[0,1]$ we have
\begin{align*}
(X_{12}f)(-1) & =  -b_{12} \,
   \Upsilon^{-1} \
 \left| \begin{matrix}
  |\alpha_{m1}'| & |\alpha_{m2}'| \\[8pt]
 |\beta_{m1}'|\, \rho_{m1}(0)\,f(1) & |\beta_{m2}'|\, \rho_{m2}(0)\, f(1)
 \end{matrix} \right|
  = -b_{12} f(1).
\end{align*}
Also
\begin{align*}
(X_{12}^*f)(1) & = -\co{b}_{12} \,
  \Upsilon^{-1} \
\left| \begin{matrix}
 |\alpha_{m1}'| & |\alpha_{m2}'| \\[8pt]
 |\alpha_{m1}'|\,f(-1) & |\alpha_{m2}'|\, f(-1)
 \end{matrix} \right|
  = 0 .
\end{align*}

\medskip

%
% case 2
%
\noindent {\bf Case} (B). \  By Theorem~\ref{tgenS} there exist
operators
\[
S_{mj} : L_{2,|r|}(0,1) \rightarrow L_{2,|r|}(-1,0), \ j=1,2,
\]
which satisfy ($S$-\ref{itS1})-($S$-\ref{itS3}) in
Theorem~\ref{tgenS} with $\imath = [0,1],\, \jmath = [-1,0]$, $a =
1$ and $b = -1$. In particular, for all $f \in \cF_{\max}[0,1]$ and
$j=1,2$,
\begin{equation*} %\label{eqpS03m2}
(S_{mj}f)(-1) = |\alpha_{mj}'| \, f(1), \ \ \ \ \ (S_{mj}^*f)(1) =
|\beta_{mj}'|\, \rho_{mj}(0) \, f(-1).
\end{equation*}
To simplify the formulas we continue to use the notation
 \[
\Upsilon : = \left| \begin{matrix} |\alpha_{m1}'| & |\alpha_{m2}'|
 \\[8pt]
 |\beta_{m1}'|\, \rho_{m1}(0) & |\beta_{m2}'|\, \rho_{m2}(0)
 \end{matrix} \right|.
 \]
Define
$$
X_{12} : L_{2,|r|}(0,1) \rightarrow L_{2,|r|}(-1,0),
 $$
by
\begin{align*}
X_{12} & =  -b_{12} \,
   \Upsilon^{-1}  \
 \left| \begin{matrix} S_{m1} & S_{m2} \\[8pt]
|\beta_{m1}'|\, \rho_{m1}(0) & |\beta_{m2}'|\, \rho_{m2}(0)
 \end{matrix} \right| .
\end{align*}
Then for all $f \in \cF_{\max}[0,1]$ we have
\begin{align*}
(X_{12}f)(-1) & =  -b_{12} \,
   \Upsilon^{-1} \
   \left| \begin{matrix}
  |\alpha_{m1}'|\,f(1) & |\alpha_{m2}'|\,f(1)
  \\[8pt]
  |\beta_{m1}'|\, \rho_{m1}(0) & |\beta_{m2}'|\, \rho_{m2}(0)
  \end{matrix} \right|
  = -b_{12} \,f(1)
\end{align*}
and for all $g \in \cF_{\max}[-1,0]$ we have
\begin{align*}
(X_{12}^*g)(1) & =  -b_{12} \,
   \Upsilon^{-1} \
 \left| \begin{matrix}
 |\beta_{m1}'|\, \rho_{m1}(0)\,g(-1)
  & |\beta_{m2}'|\, \rho_{m2}(0)\, g(-1)
  \\[8pt]
 |\beta_{m1}'|\, \rho_{m1}(0) & |\beta_{m2}'|\, \rho_{m2}(0)
 \end{matrix} \right|
  = 0.
\end{align*}

Now define the opposite off diagonal corner
 $$
X_{21} : L_{2,|r|}(-1,0) \rightarrow L_{2,|r|}(0,1),
 $$
by
\begin{align*}
X_{21} & =  b_{21} \, \Upsilon^{-1} \
 \left| \begin{matrix}
  |\alpha_{m1}'| & |\alpha_{m2}'| \\[8pt]
 S_{m1}^* & S_{m2}^*
 \end{matrix} \right| .
\end{align*}
Then  for all $f \in \cF_{\max}[-1,0]$ we have
\begin{align*}
(X_{21}f)(1) & = b_{21} \,
  \Upsilon^{-1} \
 \left| \begin{matrix}
  |\alpha_{m1}'|  & |\alpha_{m2}'| \\[8pt]
 |\beta_{m1}'|\, \rho_{m1}(0)\,f(-1) & |\beta_{m2}'|\, \rho_{m2}(0)\, f(-1)
 \end{matrix} \right|
  = b_{21} \,f(-1)
\end{align*}
and for all $g \in \cF_{\max}[0,1]$ we have
\begin{align*}
(X_{21}^*g)(-1) & =  \co{b}_{21} \,
   \Upsilon^{-1} \
 \left| \begin{matrix}
  |\alpha_{m1}'| & |\alpha_{m2}'| \\[8pt]
   |\alpha_{m1}'|\,g(1) & |\alpha_{m2}'|\, g(1) \\
 \end{matrix} \right|
  = 0.
\end{align*}

\medskip
%
% case 3
%

\noindent {\bf Case} (C). \  By Theorem~\ref{tgenS} there exists an
operator
\[
S_{m1} : L_{2,\left|r\right|}(-1,0) \rightarrow L_{2,|r|}(0,1)
\]
with the properties listed in Case (A) of this proof and there
exists an operator
\[
S_{m2}: L_{2,|r|}(0,1) \rightarrow L_{2,|r|}(-1,0)
\]
with the properties listed in Case (B).

To simplify the formulas in this part of the proof we use the
notation
 \[
\Upsilon : = \left| \begin{matrix}
  |\alpha_{m1}'| & |\beta_{m2}'|\, \rho_{m2}(0) \\[8pt]
  |\beta_{m1}'|\, \rho_{m1}(0) & |\alpha_{m2}'|
  \end{matrix} \right|.
 \]
Define
 $$
 X_{12}: L_{2,|r|}(0,1) \rightarrow L_{2,|r|}(-1,0)
 $$
by
\begin{align*}
X_{12} & =  -b_{12} \,
 \Upsilon ^{-1} \
  \left| \begin{matrix}
  |\alpha_{m1}'| & |\beta_{m2}'|\, \rho_{m2}(0) \\[8pt]
 S_{m1}^* & S_{m2}
\end{matrix} \right| .
\end{align*}
Then for all $f \in \cF_{\max}[0,1]$ we have
\begin{align*}
(X_{12}f)(-1) & = -b_{12} \,
 \Upsilon ^{-1} \
  \left| \begin{matrix}
   |\alpha_{m1}'| & |\beta_{m2}'|\, \rho_{m2}(0) \\[8pt]
  |\beta_{m1}'|\, \rho_{m1}(0)\,f(1) & |\alpha_{m2}'|\,f(1)
 \end{matrix} \right|
  = -b_{12} \,f(1)
\end{align*}
and  for all $g \in \cF_{\max}[-1,0]$ we have
\begin{align*}
(X_{12}^*g)(1) & =  -\co{b}_{12} \,
 \Upsilon ^{-1} \
  \left|\begin{matrix}
  s_{m1} & \theta_{m2}(0) \\[8pt]
 s_{m1} g(-1) & \theta_{m2}(0) g(-1)
 \end{matrix} \right|
  = 0 .
\end{align*}
The other off diagonal operator
 $$
 X_{21}: L_{2,|r|}(-1,0)\rightarrow L_{2,|r|}(0,1)
 $$
is defined as:
\begin{align*}
X_{21} & =  b_{m21} \,
  \Upsilon ^{-1} \
  \left| \begin{matrix} S_{m1} & S_{m2}^* \\[8pt]
 |\beta_{m1}'|\, \rho_{m1}(0) & |\alpha_{m2}'|
 \end{matrix} \right| .
 \end{align*}
Then  for all $f \in \cF_{\max}[-1,0]$ we have
\begin{align*}
(X_{21}f)(1) & = b_{21} \,
 \Upsilon ^{-1} \
  \left| \begin{matrix}
  |\alpha_{m1}'|\,f(-1) & |\beta_{m2}'|\, \rho_{m2}(0)\,f(-1) \\[8pt]
 |\beta_{m1}'|\, \rho_{m1}(0) & |\alpha_{m2}'|
\end{matrix} \right|
  = b_{21} \,f(-1)
\end{align*}
and  for all $g \in \cF_{\max}[0,1]$ we have
\begin{align*}
(X_{21}^*g)(-1) & = \co{b}_{21} \,
 \Upsilon ^{-1} \
 \left| \begin{matrix}
 |\beta_{m1}'|\, \rho_{m1}(0)\,g(1) & |\alpha_{m2}'|\,g(1) \\[8pt]
 |\beta_{m1}'|\, \rho_{m1}(0) & |\alpha_{m2}'|
 \end{matrix} \right|
  = 0.
\end{align*}

We conclude this part of the proof by summarizing that in each of
the three cases above we have defined operators
 $$
 X_{12}: L_{2,|r|}(0,1) \rightarrow L_{2,|r|}(-1,0) \ \ \
 \text{and} \ \ \
 X_{21}: L_{2,|r|}(-1,0)\rightarrow L_{2,|r|}(0,1)
 $$
such that
\begin{alignat*}{2}
X_{12}\cF_{\max}[0,1] &\subset \cF_{\max}[-1,0], & \ \ \  \ \
 \ \ \  X_{12}^*\cF_{\max}[-1,0] &\subset \cF_{\max}[1,0],
  \\
X_{21}^*\cF_{\max}[0,1] &\subset \cF_{\max}[-1,0], & \ \ \  \ \
 \ \ \  X_{21}\cF_{\max}[-1,0] &\subset \cF_{\max}[1,0],
\end{alignat*}
and for all $f \in \cF_{\max}[0,1]$ and $g \in \cF_{\max}[-1,0]$
we have
\begin{alignat*}{2}
(X_{12}f)(-1) & = -b_{12} \,f(1), & \ \ \  \ \
 \ \ \ (X_{12}^*g)(1) & = 0, \\
 (X_{21}^*f)(-1) & = 0, & \ \ \  \ \
 \ \ \ (X_{21}g)(1) & =  b_{21} \,f(-1).
 \hspace*{1cm}
\end{alignat*}
This completes the construction of the off-diagonal entries of
$X_{s1}$.

\medskip

\noindent{\bf 2.} \ To construct the diagonal entries we need two
self-adjoint operators $P_{1,-}$ and $P_{1,+}$ defined as follows.
Let $\phi_1:[-1,1] \rightarrow [0,1]$ be an even function with
$\phi_1\in C^1[-1,1]$ and such that
\begin{equation*} % \label{eqphi1}
 \phi_1(-1) = 1, \ \ \ \ \ \ \
 \phi_1(x) = 0 \ \ \ \text{for} \ \ \  0 \leq |x| \leq 1/2,
 \ \ \ \ \ \ \  \phi_1(1) = 1.
\end{equation*}
We now define
\begin{equation*}
P_{1,-}: L_{2,|r|}(-1,0) \rightarrow L_{2,|r|}(-1,0) \ \ \
\text{and} \ \ \  P_{1,+}: L_{2,|r|}(0,1) \rightarrow L_{2,|r|}(0,1)
\end{equation*}
by
\begin{alignat}{2} \label{eqdP1-}
(P_{1,-} f)(x) &=  f(x)\, \phi_1(x), \ \ \ f \in L_{2,|r|}(-1,0),
\ \ \ & & x \in [-1,0], \\
\intertext{and}  \label{eqdP1+}
 (P_{1,+} f)(x) &= f(x)\, \phi_1(x), \ \ \ f \in L_{2,|r|}(0,1),
 \ \ \ \ & & x \in [0,1].
\end{alignat}
These operators enjoy the following properties:
\begin{alignat*}{2} %%%\label{eqpP11-}
(P_{1,-}f)(x) & = 0, \ \ \ \ f \in L_{2,|r|}(-1,0), \ \ \ \
 &  -\tfrac{1}{2} \leq & \ x \leq 0, \\ %\label{eqpP11+}
 (P_{1,+}f)(x) &= 0, \ \ \ \ f \in L_{2,|r|}(0,1), \ \ \ \
  & 0 \leq & \ x \leq \tfrac{1}{2},
\end{alignat*}
\begin{equation*} %\label{eqpP12}
P_{1,-}\cF_{\max}[-1,0]  \subset \cF_{\max}[-1,0], \ \ \ \
 P_{1,+}\cF_{\max}[0,1]  \subset \cF_{\max}[0,1],
\end{equation*}
and
\begin{alignat*}{2}
% \label{eqpP13-}
(P_{1,-}f)(-1) &= f(-1), \ \ \ \ & f  & \in \cF_{\max}[-1,0], \\
 % \label{eqpP13+}
(P_{1,+}f)(1) & = f(1), \ \ \ \ & f  & \in  \cF_{\max}[0,1].
\end{alignat*}

Now we use Condition~\ref{cat-1} to construct the operator $X_{11}$.
As in Proposition~\ref{pwat-1}, Theorem~\ref{tgenS} implies that
there exists an operator $S_{-1}: L_{2,|r|}(-1,0) \rightarrow
L_{2,|r|}(-1,0)$ with the properties listed there. In particular for
all $f \in \cF_{\max}[-1,0]$ we have
 \begin{align*}
(S_{-1}f)(-1) & = |\alpha_{-1}'| \,f(-1), & (S_{-1}^*f)(-1) & =
|\beta_{-1}'|\, \rho_{-1}(0)\, f(-1).
\end{align*}
Since $|\alpha_{-1}'|\neq |\beta_{-1}'|\, \rho_{-1}(0)$ we can
choose complex numbers $\gamma_1$ and $\gamma_2$ such that
\[
\gamma_1 |\alpha_{-1}'| + \gamma_2 = -b_{11}-1, \ \ \ \ \ \
\co{\gamma}_1 |\beta_{-1}'|\, \rho_{-1}(0) + \co{\gamma}_2 = 1.
\]
Let $P_{1,-}$ be the operator defined in \eqref{eqdP1-}. Put
\begin{align*}
X_{11} &=  \gamma_1\ S_{-1} + \gamma_2 \ P_{1,-}.
\end{align*}
Then for all $f \in \cF_{\max}[-1,0]$ we have
\begin{equation*}
(X_{11}f)(-1)  = (-b_{11}-1) \,f(-1),  \ \ \ \ \ \ \ (X_{11}^*f)(-1)
= f(-1).
\end{equation*}
Note also that
\begin{equation*}
X_{11}\cF_{\max}[-1,0] \subset \cF_{\max}[-1,0]\ \ \ \text{and} \
\ \ X_{11}^*\cF_{\max}[-1,0] \subset \cF_{\max}[-1,0].
\end{equation*}

To construct $X_{22}$ we use Condition~\ref{cat1}. By
Theorem~\ref{tgenS} there exists a bounded operator
 \[
 S_{+1}: L_{2,|r|}(0,1) \rightarrow L_{2,|r|}(0,1)
 \]
such that
 \[
S_{+1} \cF_{\max}[0,1] \subset \cF_{\max}[0,1], \ \ \ \ S_{+1}^*
\cF_{\max}[0,1] \subset \cF_{\max}[0,1],
 \]
and for all $f \in \cF_{\max}[0,1]$,
 \begin{equation*}
(S_{+1}f)(1)  = |\alpha_{+1}'| \,f(1), \ \ \ \ \ \ \ \
(S_{+1}^*f)(1) = |\beta_{+1}'|\, \rho_{+1}(0) f(-1).
\end{equation*}
Since $|\alpha_{+1}'|\neq |\beta_{+1}'|\, \rho_{+1}(0)$ we can
choose complex numbers $\delta_1$ and $\delta_2$ such that
\[
\delta_1 |\alpha_{+1}'| + \delta_2 = -b_{11}-1, \ \ \ \ \ \
\co{\delta}_1 |\beta_{+1}'|\, \rho_{+1}(0) + \co{\delta}_2 = 1.
\]
Let $P_{1,+}$ be the operator defined in \eqref{eqdP1+}.  Put
\begin{equation*}
X_{22} =  \delta_1\, S_{+1} + \delta_2 \, P_{1,+}.
\end{equation*}
Then for all $f \in \cF_{\max}[0,1]$ we have
\begin{equation*}
(X_{22}f)(1) = (b_{22}-1) \,f(1) \ \ \ \text{and} \ \ \
(X_{22}^*f)(1) = f(1).
\end{equation*}
Note also that
\begin{equation*}
X_{22}\cF_{\max}[0,1] \subset \cF_{\max}[0,1] \ \ \ \text{and} \ \
\ X_{22}^*\cF_{\max}[0,1] \subset \cF_{\max}[0,1].
\end{equation*}

\medskip

\noindent{\bf 3.} \ Now we formally define $W_{s1} := J_0
(X_{s1}^*X_{s1} +I)$ where
\[
X_{s1} = \begin{bmatrix} X_{11} & X_{12} \\[3pt] X_{21} & X_{22}
\end{bmatrix}.
\]
To complete the proof, we verify the properties of $W_{s1}$ stated
in the theorem.  Indeed, (a) and (b) are immediate, and since
$(X_{ij}f)(x) = 0$ whenever $-\frac{1}{2} \leq x \leq \frac{1}{2}$,
(c) follows. Moreover,  each of the operators $X_{ij}$ maps
$\cF_{\max}[-1,0]$ or $\cF_{\max}[0,1]$ to $\cF_{\max}[-1,0]$ or
$\cF_{\max}[0,1]$ according  to its position in the matrix, so (d)
holds.

Finally, we check the effect of the individual components at the
boundary points $-1$ and $1$.  Evidently
\[
 X_{s1}\cF_{\max} \subset \cF_{\max}, \ \ \ \
X_{s1}^* \cF_{\max} \subset \cF_{\max}.
\]
Moreover for $f,g \in \cF_{\max}$ we have
\begin{align*}
\begin{bmatrix}
(X_{s1}f)(-1) \\[5pt] (X_{s1}f)(1)
\end{bmatrix} & =
\begin{bmatrix}
(X_{11}f)(-1) + (X_{12}f)(-1) \\[5pt] (X_{21}f)(1) +  (X_{22}f)(1)
\end{bmatrix} = \begin{bmatrix}
(-b_{11}-1)f(-1) - b_{12}f(1) \\[5pt] b_{21}f(-1) +
(b_{22}-1)f(1)
\end{bmatrix}
\end{align*}
and
\begin{align*}
\begin{bmatrix}
(X_{s1}^*g)(-1) \\[5pt] (X_{s1}^*g)(1)
\end{bmatrix} & =
\begin{bmatrix}
(X_{11}^*g)(-1) + (X_{21}^*g)(-1) \\[5pt] (X_{12}^*g)(1) +
(X_{22}^*g)(1)
\end{bmatrix} = \begin{bmatrix}
g(-1) + 0 \\[5pt] 0 + g(1)
\end{bmatrix}.
\end{align*}
Substituting $g= X_{s1}f \in \cF_{\max}$, we get
\begin{align*}
\begin{bmatrix}
(X_{s1}^*X_{s1}f)(-1) \\[5pt] (X_{s1}^* X_{s1}f)(1)
\end{bmatrix} & =
\begin{bmatrix}
(-b_{11}-1)f(-1) - b_{12}f(1) \\[5pt] b_{21}f(-1) +
(b_{22}-1)f(1)
\end{bmatrix}.
\end{align*}
With $Y_{s1} = X_{s1}^*X_{s1} +I$ we have
\begin{align*}
\begin{bmatrix}
(Y_{s1}f)(-1) \\[5pt] (Y_{s1}f)(1)
\end{bmatrix} & =
\begin{bmatrix}
-b_{11}f(-1) - b_{12}f(1) \\[5pt] b_{21}f(-1) + b_{22}f(1)
\end{bmatrix} = \begin{bmatrix}
-b_{11} & -b_{12} \\[5pt] b_{21} & b_{22}
\end{bmatrix} \begin{bmatrix}
 f(-1) \\[5pt] f(1) \end{bmatrix},
\end{align*}
which proves (e) since $W_{s1} = J_0Y_{s1}$.
\end{proof}

\begin{remark}
Notice that the operators $W_{-1}$ and $W_{+1}$ from
Propositions~\ref{pwat-1} and \ref{pwat1} satisfy
\[
\begin{bmatrix}
(W_{-1}f)(-1) \\[5pt]  \rule{0pt}{0pt} (W_{-1}f)(1)
\end{bmatrix} =  {\begin{bmatrix}
b & 0 \\[5pt]
0 & 1
\end{bmatrix}} \begin{bmatrix}
f(-1) \\[5pt]  \rule{0pt}{0pt} f(1)
\end{bmatrix} \ \ \ \text{and} \ \ \
\begin{bmatrix}
(W_{+1}f)(-1) \\[5pt]  \rule{0pt}{0pt} (W_{+1}f)(1)
\end{bmatrix} =  {\begin{bmatrix}
-1 & 0 \\[5pt]
0 & b
\end{bmatrix}} \begin{bmatrix}
f(-1) \\[5pt]  \rule{0pt}{0pt} f(1)
\end{bmatrix},
\]
respectively, with arbitrary $b \in \nC$.  A stronger conclusion is
contained in Theorem~\ref{twac1}~(\ref{itwac1e}) under stronger
assumptions.
\end{remark}

\section{Two essential or two non-essential boundary conditions}
  \label{sk=0}

The first theorem of this section deals with the case of two
non-essential boundary conditions.

\begin{theorem} \label{tk0sebc}
Assume that the following two conditions are satisfied.
\begin{enumerate}[{\rm (a)}]
\item
  $\M{N}_n = \begin{bmatrix}
  1 & 0 \\ 0 & 1   \end{bmatrix}$.
\item
The coefficients $p$ and $r$ satisfy Condition~{\rm~\ref{c0}}.
\end{enumerate}
Then there is a basis for each root subspace of $A$, so that the
union of all these bases is a Riesz basis of $L_{2,|r|} \oplus
\nC_{|\Delta|}$.
\end{theorem}

\begin{proof}
By \eqref{fd1}, the form domain of $\w{A}$ is given as
\begin{align*}
\fdom(\w{A}) & = \left\{ \begin{bmatrix} f \\[5pt] \V{v}
\end{bmatrix} \in
\begin{matrix}
L_{2,r} \\ \oplus \\ \nC^2_{\Delta} \end{matrix} \, : \, f \in
\cF_{\max}, \, \V{v} \in \nC^2 \right\}.
\end{align*}
Recalling $W_0$ from Theorem~\ref{twat0}, we easily see that the
operator
\[
\w{W} = \begin{bmatrix} W_0 \ & 0 \\[5pt]  0 \ & \M{\Delta}^{\!\!^{-1}}
\end{bmatrix} \ : \ \begin{matrix}L_{2,r} \\ \oplus  \\
\nC_{\M{\Delta}}^2
\end{matrix}
\ \ \rightarrow \ \ \begin{matrix}L_{2,r} \\ \oplus  \\
\nC_{\M{\Delta}}^2 \end{matrix}.
\]
is bounded, boundedly invertible and positive in the Krein space
$L_{2,r}\oplus \nC^2_{\Delta}$. A simple verification shows that
$\w{W} \fdom(\w{A}) \subset \fdom(\w{A})$ so the theorem follows
from Theorem~\ref{W}.
\end{proof}

We now consider the case of two essential conditions.

\begin{theorem} \label{tn=k2}
Assume that the following three conditions are satisfied.
\begin{enumerate}[{\rm (a)}]
\item \label{tn=k21}
 $\M{N}_e =  \begin{bmatrix} 1 & 0  \\[2pt]
 0 & 1 \end{bmatrix}$.
\item  \label{tn=k22}
The coefficients $p$ and $r$ satisfy Condition~{\rm~\ref{c0}}.
\item  \label{tn=k23}
One of the following holds:
\begin{enumerate}[{\rm (i)}]
\item
$\M{\Delta} > 0$ and the coefficients $p$ and $r$ satisfy
Condition~{\rm~\ref{cat-1}}.
\item
$\M{\Delta} < 0$ and the coefficients $p$ and $r$ satisfy
Condition~{\rm~\ref{cat1}}.
\item
the coefficients $p$ and $r$ satisfy Conditions~{\rm \ref{cat-1}},
{\rm \ref{cat1}} and~{\rm~\ref{mc1}}.
\end{enumerate}
\end{enumerate}
Then there is a basis for each root subspace of $A$, so that the
union of all these bases is a Riesz basis of $L_{2,|r|} \oplus
\nC_{|\Delta|}$.
\end{theorem}

\begin{proof}
Define the following two
Krein spaces:
\begin{equation*} % \label{eqK0K1}
 \cK_0  := L_{2,r}\!\left(-\tfrac{1}{2},\tfrac{1}{2}\right),  \ \ \
 \ \ \   \cK_1 := L_{2,r}(-1,-\tfrac{1}{2})[\dot{+}]
 L_{2,r}(\tfrac{1}{2},1).
\end{equation*}
Extending functions in $\cK_0$ and $\cK_1$ by zero, we consider
the spaces $\cK_0$ and $\cK_1$ as subspaces of $L_{2,r}$. Then
 \[
L_{2,r} = \cK_0[\dot{+}]\cK_1.
 \]

As in the previous proof our goal is to construct $\w{W}: L_{2,r}
\oplus \nC^2_{\M{\Delta}} \to L_{2,r} \oplus \nC^2_{\M{\Delta}}$.
The first step is to define $W_{01}:L_{2,r}\to L_{2,r}$. We proceed
by considering each case in (c)  separately.

\medskip

\noindent (i) Let $W_0$ be the operator constructed in
Theorem~\ref{twat0} and let $W_{-1}$ be the operator constructed in
Proposition~\ref{pwat-1} with $b = 1$.  Property (\ref{ipwat-1c}) in
Theorem~\ref{twat0} and Proposition~\ref{pwat-1} imply that $\cK_0$
and $\cK_1$ are invariant under $W_0$ and $W_{-1}$.  Since we chose
$b = 1$, we have $(W_{-1}f)(-1) = f(-1)$ and $(W_{-1}f)(1) = f(1)$.
Define
\begin{equation} \label{eqdW011}
W_{01} : = W_0|_{_{\cK_0}}[\dot{+}]W_{-1}|_{_{\cK_1}}.
\end{equation}
Since $W_0$ and $W_{-1}$ are bounded, boundedly invertible and
positive in the Krein space $L_{2,r}$, so is the the operator
$W_{01}$.
% is bounded, boundedly invertible and
% positive in the Krein space $L_{2,r}$.
Also, $W_{01} \cF_{\max} \subset \cF_{\max}$ and
\begin{equation} \label{eqW01-}
\begin{bmatrix} (W_{01}f)(-1) \\[5pt] (W_{01}f)(1)
\end{bmatrix} = \begin{bmatrix} f(-1)
\\[5pt] f(1) \end{bmatrix}.
\end{equation}

\medskip

\noindent (ii) Instead of $W_{-1}$ in (i), we use the operator
$W_{+1}$ constructed in Proposition~\ref{pwat1} with $b = -1$.
Redefining the operator $W_{01}$ as
\begin{equation} \label{eqdW012}
W_{01} : = W_0|_{_{\cK_0}}[\dot{+}]W_{+1}|_{_{\cK_1}}.
\end{equation}
we see that it is again bounded, boundedly invertible, and positive
in the Krein space $L_{2,r}, \; W_{01} \cF_{\max} \subset
\cF_{\max}$ and (since we use $b = -1$)
\begin{equation} \label{eqW01+}
\begin{bmatrix} (W_{01}f)(-1) \\[5pt] (W_{01}f)(1)
\end{bmatrix} = - \begin{bmatrix} f(-1)
\\[5pt] f(1) \end{bmatrix}.
\end{equation}

\medskip

\noindent (iii) This time we replace  $W_{-1}$ from (i) by $W_{s1}$
from Theorem~\ref{twac1}, so we define the operator
\begin{equation}\label{W0s1}
W_{01} : = W_0|_{_{\cK_0}}[\dot{+}]W_{s1}|_{_{\cK_1}},
\end{equation}
which is again
bounded, boundedly invertible and positive in the Krein space
$L_{2,r}$.  Also, $W_{01} \cF_{\max} \subset \cF_{\max}$ and
\begin{equation} \label{eqW01D}
\begin{bmatrix}
(W_{s1}f)(-1) \\[5pt]  \rule{0pt}{0pt} (W_{s1}f)(1)
\end{bmatrix} =  \M{\Delta}^{-1} \begin{bmatrix}
f(-1) \\[5pt]  \rule{0pt}{0pt} f(1)
\end{bmatrix}.
\end{equation}

Finally we define $\w{W}: L_{2,r} \oplus
\nC^2_{\M{\Delta}} \to L_{2,r} \oplus \nC^2_{\M{\Delta}}$ by
\begin{equation}\label{W-}
\w{W} = \begin{bmatrix} W_{01} \ & 0 \\[5pt]  0 \ & I
\end{bmatrix}
\end{equation}
in case (c)(i),
\begin{equation}\label{W+}
\w{W} = \begin{bmatrix} W_{01} \ & 0 \\[5pt]  0 \ & -I
\end{bmatrix}
\end{equation}
in case (c)(ii), and
\begin{equation}\label{WD}
\w{W} = \begin{bmatrix} \ W_{01} \ & 0 \  \\[5pt]
\rule{0pt}{0pt} \ 0 \ & \M{\Delta}^{-1} \
\end{bmatrix}
\end{equation}
in case (c)(iii).

By \eqref{fd3}, the form domain of $\w{A}$ is
\begin{equation*}
 \fdom(A) = \left\{ \begin{bmatrix} f \\[5pt] f(-1) \\ f(1)
\end{bmatrix} \in
\begin{matrix} L_{2,r} \\ \oplus \\ \nC_{\Delta}^2 \end{matrix}
 \ : \ f \in \cF_{\max} \right\}.
\end{equation*}
A straightforward verification shows that in each case \eqref{W-},
\eqref{W+}, and \eqref{WD}, $\w{W}$ is a bounded, boundedly
invertible, positive operator in the Krein space $L_{2,r}\oplus
\nC^2_{\Delta}$. Moreover $\w{W}\fdom(\w{A}) \subset \fdom(\w{A})$
via \eqref{eqW01-}, \eqref{eqW01+} or \eqref{eqW01D}. Now the
theorem follows from Theorem~\ref{W}.
\end{proof}

\begin{example} \label{esbe}
Consider the eigenvalue problem
\begin{align*}
-f'' & = \lambda \, r \, f \\
f'(1) & = \lambda f(-1) \\
-f'(-1) & = \lambda f(1),
\end{align*}
where $r(x) = \sgn x, x \in [-1,1],$ as in our example in the
Introduction.  Then clearly
 $\M{N}_e =  \begin{bmatrix} 1 & 0  \\[2pt]
 0 & 1 \end{bmatrix}$,
giving (a) in Theorem~\ref{tn=k2} and (b) follows from the note after
Condition~\ref{c0}. Moreover, an easy computation gives
 $\Delta = \begin{bmatrix} 0 & 1 \\
           1 & 0 \end{bmatrix}$,
which is indefinite. Condition (c) now follows from
Examples~\ref{e42} and~\ref{r42}, so Theorem~\ref{tn=k2} applies.

On the other hand, if instead we take $r$ as in Example~\ref{e44},
then as we have seen, Condition~\ref{mc1} fails and hence so does
(c)(iii) in Theorem~\ref{tn=k2}. Therefore Theorem~\ref{tn=k2} gives
no conclusion about a Riesz basis for this amended case.
\end{example}

\section{%
 One essential and one non-essential boundary condition}
 \label{sn2k1}

The main result of this section is the following theorem. Its proof
will occupy the most of the section and then we will proceed to some
examples.

\begin{theorem} \label{tn2k1}
Assume that the following three conditions are satisfied.
\begin{enumerate}[{\rm (a)}]
\item \label{tn2k11}
$ \displaystyle  \M{N} =  \begin{bmatrix} u & v & 0 & 0 \\[3pt]
* & * & * & 1
\end{bmatrix} \ \ \text{or} \ \
\M{N} =  \begin{bmatrix} u & v & 0 & 0 \\[3pt]
* & * & 1 & 0
\end{bmatrix}$, where $|u|^2 + |v|^2 > 0$ and the asterisks stand
for arbitrary complex numbers.
\item
The coefficients $p$ and $r$ satisfy Condition~{\rm~\ref{c0}}.
\item
One of the following holds.
\begin{enumerate}[{\rm (i)}]
\item
$u = 1, v = 0$ and the coefficients $p$ and $r$ satisfy
Condition~{\rm~\ref{cat-1}}.
\item
$u = 0, v = 1$ and the coefficients $p$ and $r$ satisfy
Condition~{\rm~\ref{cat1}}.
\item
$uv \neq 0$ and the coefficients $p$ and $r$ satisfy
Conditions~{\rm~\ref{cat-1}}, {\rm \ref{cat1}} and~{\rm~\ref{mc1}}.
\end{enumerate}
\end{enumerate}
Then there is a basis for each root subspace of $A$, so that the
union of all these bases is a Riesz basis of $L_{2,|r|} \oplus
\nC_{|\Delta|}$.
\end{theorem}

\begin{proof}
It follows from (\ref{tn2k11}) that the form domain of $\w{A}$ is
\begin{align}
\fdom(\w{A})  & = \left\{ \begin{bmatrix} f \ \ \  \ \  \   \\[2pt]
\rule{0pt}{0pt} u f(-1) + v f(1) \\[2pt]
\rule{0pt}{0pt} z \ \ \  \ \ \
\end{bmatrix} \ \in \
\begin{matrix}L_{2,r} \\ \oplus  \\ \nC^2_{\M{\Delta}} \end{matrix}
 \ : \  f \in \cF_{\max}, \ z \in \nC
\right\}. \label{fd}
\end{align}
It is no restriction if we scale the first boundary condition so
that
\begin{equation} \label{equv1}
|u|^2 + |v|^2 =1.
\end{equation}

As in the previous proofs we shall construct $\w{W}:L_{2,r}\oplus
\nC^2_{\M{\Delta}} \to L_{2,r}\oplus \nC^2_{\M{\Delta}}$ in blocks.
We divide the proof into three parts and two lemmas.

\medskip

\noindent{\bf 1.} \ First we define a bounded operator
$W_{01}:L_{2,r} \to L_{2,r}$ such that
\begin{gather} \label{eqpW1}
 J_0 W_{01} > I,  \\ \label{eqpW2}
W_{01} \cF_{\max} \subset \cF_{\max}, \\
  \label{eqpW}
u(W_{01}f)(-1) + v(W_{01}f)(1) = 0, \ \ \ \ f \in \cF_{\max}.
\end{gather}

We distinguish the three cases in (c) above.

\medskip

\noindent (i) \ As in the proof of Theorem~\ref{tn=k2}(i), we define
$W_{01}$ by \eqref{eqdW011}, but now using $b = 0$ instead of $b =
1$. Then $W_{01}$ is a bounded operator in the Krein space
$L_{2,r}$, and it satisfies \eqref{eqpW2} and $(W_{01}f)(-1) = 0, \,
(W_{01}f)(1)= f(1)$ and hence \eqref{eqpW}. Inequality \eqref{eqpW1}
follows from \eqref{eqdW011}, Theorem~\ref{twat0}(\ref{itwat0b}) and
Proposition~\ref{pwat-1}(\ref{ipwat-1b}).

\medskip

\noindent (ii) \ This time we define $W_{01}$ by \eqref{eqdW012},
but now using $b = 0$ instead of $b = -1$. Then $W_{01}$ is a
bounded operator in the Krein space $L_{2,r}$, it satisfies
\eqref{eqpW2} and $(W_{01}f)(-1) = -f(-1), \, (W_{01}f)(1)= 0$ and
hence \eqref{eqpW}. In this case inequality \eqref{eqpW1} follows
from \eqref{eqdW012}, Theorem~\ref{twat0}(\ref{itwat0b}) and
Proposition~\ref{pwat1}(\ref{ipwat1b}).

\medskip

\noindent (iii) \ We now define $W_{01}$ as in the proof of
Theorem~\ref{tn=k2}(iii), but instead of using $\M{\Delta}^{-1}$ in
\ref{eqW01D} we use the zero $2\times 2$ matrix $\M{0}$.  Then
$W_{01}$ is a bounded operator in the Krein space $L_{2,r}$, it
satisfies \eqref{eqpW2} and $(W_{01}f)(-1) = 0, \, (W_{01}f)(1)= 0$
and hence \eqref{eqpW}. Inequality \eqref{eqpW1} follows from
\eqref{W0s1}, Theorem~\ref{twat0} (\ref{itwat0b}) and
Theorem~\ref{twac1} (\ref{itwac1b}).

\medskip

\noindent{\bf 2.} \ Next we define an integral operator $K$ which
will be a perturbation of $W_{01}$.

\noindent{\bf 2.1.} \ We start by writing the inverse of the matrix
$\M{\Delta}$ in the form
\[
\M{\Delta}^{-1} =  \begin{bmatrix} \eta_{11} & \eta_{12} \\[3pt]
\co{\eta}_{12} & \eta_{22}
\end{bmatrix},
\]
and setting $\eta : = \max\{|\eta_{11}|,|\eta_{12}| \} > 0$, with
$\delta_2 \geq \delta_1 > 0$ as the eigenvalues of $|\M{\Delta}|$.
We also define three positive constants
\begin{align} \nonumber % \label{eq3in1}
\alpha & := \frac{ \delta_2 }{1+2\|r\|_1\delta_2 \eta^2} , \\
 \label{eq3in2}
c & := \frac{\alpha}{2 \delta_2} \ \sqrt{\frac{\delta_1}{2}},   \\
 \label{eq3in3}
\kappa  & :=  \frac{ 2 \delta_2 \eta^2 \|r\|_1 }{1+2\|r\|_1\delta_2
\eta^2} = 2 \, \alpha \eta^2 \|r\|_1.
\end{align}
Notice that
\begin{equation} \label{eqidn}
1 - \kappa = \dfrac{1}{1+2\|r\|_1\delta_2 \eta^2} =
\dfrac{\alpha}{\delta_2}.
\end{equation}

\bigskip

\noindent{\bf 2.2.} Since $r$ is integrable over $[-1,1]$, we there
exists $\gamma \in [0,1)$ such that
\begin{equation} \label{eqidn1}
- \int_{-1}^{-\gamma} r + \int_\gamma^1 r \leq \left(\frac{c}{\alpha
\eta} \right)^2.
\end{equation}
Noting that $p^{-1/2}\in L_2(0,1) \subset L_1(0,1)$ we can define
\[
\phi(x) = \int_0^x p^{-1/2}\chi_{[\gamma,1]}, \ \ \ x\in [0,1].
\]
Extending $\phi$ as an even function over $[-1,1]$ we see that $\phi
\in \cF_{\max}$. Since $\phi(1)$ is a positive real number, we
define $\psi = \phi/\phi(1)$. Clearly $\psi:[-1,1] \to [0,1]$ is an
even function in $\cF_{\max}$ such that
\begin{equation} \label{eqpsib}
\psi(-1) = 1, \ \ \ \ \psi(0) = 0, \ \ \ \ \psi(1) = 1,
\end{equation}
and, by \eqref{eqidn1},
\begin{equation} \label{eqnpsi}
\|\psi\|_{2,|r|} \leq \frac{c}{\alpha \eta}.
\end{equation}

\noindent{\bf 2.3.} Define
\begin{equation} \label{eqdpsi}
\psi_j(x) = \begin{cases}
 \alpha\, \eta_{1j}\, \co{u}\ \psi(x),  \ \ \ & x \in [-1,0), \\[5pt]
\alpha\, \eta_{1j}\, \co{v}\ \psi(x),  \ \ \ & x \in [0,1].
 \end{cases}
\end{equation}
Since $\psi \in \cF_{\max}$ and $\psi(0) = 0$, the functions
$\psi_1$ and $\psi_2$ belong to $\cF_{\max}$.  Set
\[
\omega(x) := \eta_{11}\,
\co{\psi_1(x)} + \eta_{12}\, \co{\psi_2(x)}, \ \ \ x \in [-1,1],
\]
and define $k: [-1,1] \times [-1,1] \rightarrow \nC$ by
\begin{equation} \label{eqdk}
k(x,t) = \begin{cases}
 u \, \co{\omega(x)} \ \ \ & \text{if} \ \ \ t \leq -|x|, \\
\co{v} \, \omega(t) \ \ \ & \text{if} \ \ \ x > |t|, \\
v \, \co{\omega(x)} \ \ \ & \text{if} \ \ \ t \geq |x|, \\
 \co{u} \, \omega(t)\ \ \ & \text{if} \ \ \ x < -|t|.  \\
\end{cases}
\end{equation}
By the definitions of $\psi_1, \psi_2$ and $\omega$, since $\psi$ is
a nonnegative even function, for all $x \in [0,1]$ we have
 \begin{equation} \label{eq3om}
\co{u}\, \omega(-x), \ \co{v}\, \omega(x) \in \nR, \ \ \
\text{and} \ \ \ \co{v}\, \omega(-x) = u \, \co{\omega(x)}.
 \end{equation}
Since $\omega$ is continuous, it follows from \eqref{eq3om} and
\eqref{eqdk} that $k$ is a continuous function.  Moreover, by
\eqref{equv1} and \eqref{eqdpsi},
\[
 |\omega(t)| < \eta \, \eta \alpha +
 \eta \, \eta \alpha  = 2\eta^2 \alpha.
\]
Therefore \eqref{eq3in3} shows that
\begin{equation} \label{equbk}
|k(x,t)| \leq 2\eta^2 \alpha =
\frac{\kappa}{\|r\|_1} .
\end{equation}

The first of our two lemmas is as follows.

\begin{lemma}\label{K}
Let $K: L_{2,r} \rightarrow L_{2,r}$ be the integral operator
defined by
\[
(K f)(x) := \int_{-1}^{1} k(x,t)\, f(t) \, r(t) \, dt, \ \ \ f\in
L_{2,r}.
\]
Then
\begin{enumerate}[\rm{(}I\rm{)}]
 \item
 The operator $K$ is bounded and self-adjoint on $L_{2,r}$ and
$\|K\|_{2,|r|} \leq \kappa$.
 \item
The range of $K$ is contained in $\cF_{\max}$.
\end{enumerate}
\end{lemma}

\begin{proof}
(I) \ We first note that for $f$ in $L_{2,r}$ the function $fr$ is
integrable on $(-1,1)$. In fact
\begin{equation} \label{eqfr}
\begin{split}
 \int_{-1}^{1} |fr| & = \int_{-1}^{1} |r|^{1/2}
\bigl(|f||r|^{1/2}\bigr) \\
 & \leq \left( \int_{-1}^{1} |r| \right)^{1/2}
 \left( \int_{-1}^{1} |f|^{2} |r| \right)^{1/2} = \| r \|_1^{1/2}
 \|f\|_{2,r}.
\end{split}
\end{equation}
For $f \in L_{2,|r|}$ we calculate
\begin{align*}
\| Kf \|_{2,|r|}^2 & \leq  \int_{-1}^{1} \int_{-1}^{1} |k(x,t)|\,
|f(t)| \, |r(t)| dt \, \int_{-1}^{1} |k(x,s)|\,
|f(s)| \, |r(s)|\,  ds  \, |r(x)|dx \\
 & \leq \frac{\kappa^2}{\|r\|_1^2} \, \int_{-1}^{1}
  \left( \int_{-1}^{1} |f| \, |r|
 \right)^2 |r(x)| dx \leq  \kappa \|f\|_{2,|r|}^2,
\end{align*}
by virtue of \eqref{equbk} and \eqref{eqfr}.  Thus $\|K\|_{2,|r|}
\leq \kappa$, so $K$ is bounded, and  self-adjointness
follows from \eqref{eqdk} since $k(x,t) =
\co{k(t,x)}, \, x,t \in [-1,1]$.

\medskip

\noindent (II) \ Let $f \in L_{2,r}$.  By definition, for $-1 \leq x
< 0$,
\begin{multline*} %\label{eqK<0}
(Kf)(x) =
  u\, \co{\omega(x)}\, \int_{-1}^{x} (f r)(t)dt
  +   \co{u}  \int_{x}^{-x} \bigl(\omega f r\bigr)(t)dt
 +  v\, \co{\omega(x)}\, \int_{-x}^{1} (f r )(t)dt
\end{multline*}
 and, for $0 < x \leq 1$,
\begin{multline*} %\label{eqK>0}
 (Kf)(x)  =
 u\, \co{\omega(x)}\, \int_{-1}^{-x} (f r)(t)dt
  +   \co{v}  \int_{-x}^{x}\bigl(\omega f r\bigr)(t)dt
  + v\, \co{\omega(x)}\, \int_{x}^{1} (f r )(t)dt  .
\end{multline*}
The function $fr$ is integrable on $(-1,1)$ by \eqref{eqfr}. Since
$\omega \in \cF_{\max}$ the function $\omega f r$ is also
integrable on $(-1,1)$. Moreover,
\[
\lim\limits_{x\uparrow 0} (Kf)(x) = \lim\limits_{x\downarrow 0}
(Kf)(x) = (Kf)(0) = 0.
\]
Therefore for each $f \in L_{2,|r|}$ the function $Kf$ is
absolutely continuous on $[-1,1]$. For almost all $x \in [-1,0)$,
we have
\begin{multline*}
(Kf)'(x)  =
  u\, \co{\omega}'(x)\, \int_{-1}^{x} (f r)(t)dt
  +  v\, \co{\omega}'(x)\, \int_{-x}^{1} (f r )(t)dt \\
   +  u\, \co{\omega(x)}\, (f r)(x)
 -  \co{u}\, \omega(x)\, \bigl(f r\bigr)(x)  -
  \co{u} \, \omega(-x) \, \bigl(f r\bigr)(-x)
 + v\, \co{\omega(x)}\,(f r )(-x),
\end{multline*}
and, for almost all $x \in (0,1]$,
\begin{multline*}
(Kf)'(x)  = u \co{\omega}'(x) \int_{-1}^{-x} (f r)(t)dt
 + v \co{\omega}'(x) \int_{x}^{1} (f r )(t)dt  \\
  -  u\, \co{\omega(x)} \, (f r)(-x) + \co{v} \, \omega(-x) \,
  \bigl( f r\bigr)(-x) + \co{v}\, \omega(x)\, \bigl(f r\bigr)(x)
  -  v\, \co{\omega(x)}\, (f r)(x) .
\end{multline*}
By \eqref{eq3om} the terms not involving integrals in the above
two equations cancel in pairs.  Thus $Kf \in \cF_{\max}$  for all
$f \in L_{2,|r|}$  since $\co{\omega} \in \cF_{\max}$. This
completes the proof of the lemma.
\end{proof}

\medskip

\noindent{\bf 3.} \ We create off-diagonal blocks for $\w{W}$
by means  of the operator $Z:
\nC_{\M{\Delta}}^2 \rightarrow L_{2,r}$ which we define by
\[
Z \V{a}:= a_1\, \psi_1 + a_2 \, \psi_2, \ \ \ \
\V{a}=\begin{bmatrix} a_1 \\[3pt] a_2
\end{bmatrix}  \in \nC^2.
\]
The adjoint $Z^{[*]}: L_{2,r} \rightarrow \nC_{\M{\Delta}}^2$ of
$Z$ is given by
\[
Z^{[*]}f = \M{\Delta}^{-1} \begin{bmatrix} [f,\psi_1] \\[3pt]
[f,\psi_2]
\end{bmatrix}, \ \ \ f \in L_{2,r}.
\]
Equalities \eqref{equv1}, \eqref{eqnpsi} and \eqref{eqdpsi} yield
$\|\psi_1\|_{2,|r|} \leq  c$ and $\|\psi_2\|_{2,|r|} \leq  c$.
 Therefore
\begin{align*}
\int_{-1}^1 \bigl| Z \V{a} \bigr|^2 |r| & \leq 2
\bigl( |a_1|^2 \| \psi_1\|^2_{2,|r|} + |a_2|^2 \| \psi_2
\|^2_{2,|r|}
\bigr) \\
& \leq 2 c^2 \V{a}^* \V{a} \leq 2 \frac{c^2}{\delta_1} \V{a}^*
 |\M{\Delta} | \V{a} .
\end{align*}
Consequently, by \eqref{eq3in2},
\begin{equation}\label{nZ}
\|Z\| = \|Z^{[*]}\| \leq c \sqrt{\frac{2}{\delta_1}} =
\frac{\alpha}{2 \delta_2}.
\end{equation}

\bigskip

The second lemma we need is as follows.

\begin{lemma}\label{fin}
Let the operator $W : L_{2,r} \oplus \nC_{\M{\Delta}}^2 \rightarrow
L_{2,r} \oplus \nC_{\M{\Delta}}^2$ be defined by
\[
\w{W} := \begin{bmatrix}   W_{01} + K & Z \\[5pt] Z^{[*]} & \alpha \,
\M{\Delta}^{-1}
\end{bmatrix}.
\]
Then
\begin{enumerate}[\rm{(}I\rm{)}]
 \item
$\w{W}$ is bounded and uniformly positive on $L_{2,r}\oplus
\nC^2_{\M{\Delta}}$.
 \item
 $\w{W} \fdom(\w{A}) \subset \fdom(\w{A})$.
\end{enumerate}
\end{lemma}
\begin{proof}
(I) The operator $\w{W}$ is bounded since each of its components is
bounded. To prove that $\w{W}$ is uniformly positive, we shall show
that the operator $\w{J}\, \w{W}$ is uniformly positive in the
Hilbert space $L_{2,|r|}\oplus \nC^2_{\M{|\Delta}|}$.  From
Lemma~\ref{K}, $\|K\| = \|JK\| \leq \kappa$ and
\[
\left\| \begin{bmatrix} 0 & Z \\ Z^{[*]} & 0
\end{bmatrix} \right\|
 = \left\| \w{J} \begin{bmatrix} 0 & Z \\ Z^{[*]} & 0
\end{bmatrix} \right\| \leq \frac{\alpha}{2 \delta_2}
\]
follows from \eqref{nZ}.  Thus
\begin{equation*}
\begin{split}
\left\la \w{J}\, \w{W} \begin{bmatrix} f \\ \V{a}
\end{bmatrix}, \begin{bmatrix} f \\ \V{a}
\end{bmatrix} \right\ra & =
 \left\la  \begin{bmatrix}
 J_0W_{01} & 0 \\ 0 & \alpha \M{|\Delta|^{-1}}
\end{bmatrix}\begin{bmatrix} f \\ \V{a} \end{bmatrix},
\begin{bmatrix} f \\ \V{a}
\end{bmatrix} \right\ra  \\  & \ \ \ \ \ \ \ \
      +  \left\la  \begin{bmatrix}
J_0K & 0 \\ 0 & 0 \end{bmatrix}\begin{bmatrix} f
\\ \V{a} \end{bmatrix}, \begin{bmatrix} f \\ \V{a}
\end{bmatrix} \right\ra +
 \left\la \w{J} \begin{bmatrix}
0 &  Z \\ Z^{[*]}  & 0 \end{bmatrix}\begin{bmatrix} f
\\ \V{a} \end{bmatrix}, \begin{bmatrix} f \\ \V{a}
\end{bmatrix} \right\ra \\
& = \la J_0 W_{01} f,f \ra + \alpha\, \V{a}^*
\V{a} + \la J_0K f,f\ra + \left\la \w{J} \begin{bmatrix} 0 &  Z \\
Z^{\la * \ra}  & 0
\end{bmatrix}\begin{bmatrix} f
\\ \V{a} \end{bmatrix}, \begin{bmatrix} f \\ \V{a}
\end{bmatrix} \right\ra \\
& \geq  \la f,f\ra + \frac{\alpha}{\delta_2} \, \V{a}^*
|\M{\Delta}| \V{a} - \kappa \la f,f\ra -  \frac{\alpha}{2\delta_2}
\,
 \bigl(\la f,f \ra + \V{a}^*|\M{\Delta}|\V{a}\bigr) \\
& \geq  \left( 1 - \kappa - \frac{\alpha}{2\delta_2} \right) \la
f,f \ra + \left( \frac{\alpha}{\delta_2} -
\frac{\alpha}{2\delta_2} \right)\,
\V{a}^*|\M{\Delta}|\V{a} \\
& = \left( \frac{\alpha}{\delta_2} - \frac{\alpha}{2\delta_2}
\right) \la f,f \ra + \frac{\alpha}{2\delta_2} \,
\V{a}^*|\M{\Delta}|\V{a}
 \text{ \hspace*{1in} (by \eqref{eqidn})} \\
& = \frac{\alpha}{2\delta_2}  \, \bigl( \la f,f \ra +
\V{a}^*|\M{\Delta}|\V{a} \bigr) \\
 & = \frac{\alpha}{2\delta_2}  \, \left\la \begin{bmatrix} f
\\ \V{a} \end{bmatrix}, \begin{bmatrix} f \\ \V{a}
\end{bmatrix} \right\ra,
\end{split}
\end{equation*}
as required.

\medskip

\noindent (II) \ We start with the identity
\begin{equation} \label{eqebcK}
u (Kf)(-1) + v (Kf)(1) = \eta_{11}\, [f,\psi_1] + \eta_{12}\,
[f,\psi_2], \ \ \ f \in L_{2,|r|},
\end{equation}
which follows from the calculation
\begin{align*}
  u (Kf)(-1) & + v (Kf)(1) \\
  & = u \int_{-1}^{1} k(-1,t)\, f(t) \, r(t) \, dt
  + v \int_{-1}^{1} k(1,t)\, f(t) \, r(t) \, dt \\
 & = u \int_{-1}^{1} \co{u}
  \bigl(\eta_{11} \co{\psi}_1(t) + \eta_{12} \co{\psi}_2(t) \bigr)\,
   f(t) \, r(t) \, dt \\
 & \ \hspace{1cm} \ \ \  + v \int_{-1}^{1}\co{v}
  \bigl(\eta_{11} \co{\psi}_1(t) + \eta_{12} \co{\psi}_2(t) \bigr)\,
   f(t) \, r(t) \, dt \\
 & = |u|^2 \eta_{11} [f,\psi_1] + |u|^2 \eta_{12} [f,\psi_2] +
 |v|^2 \eta_{11} [f,\psi_1] + |v|^2 \eta_{12} [f,\psi_2]\\
 & = \eta_{11} [f,\psi_1] + \eta_{12} [f,\psi_2].
\end{align*}

By \eqref{fd}, the general element of $\fdom(\w{A})$
takes the form
\[
\begin{bmatrix} f \ \ \  \ \  \   \\[2pt]
u f(-1) + v f(1) \\[2pt]
 z \ \ \  \ \ \
\end{bmatrix}
\]
where $f \in \cF_{\max}$ and $z \in \nC$.
Applying $\w{W}$ to this vector we obtain
 \[
w:= \begin{bmatrix}
g \\[2pt]
\eta_{11}\,[f,\psi_1]+\eta_{12}\,[f,\psi_2]+
\alpha\,\eta_{11}\,\bigl(u f(-1)
  + v f(1) \bigr) +\alpha\,\eta_{12}\, z \\[2pt]
*
\end{bmatrix},
 \]
where
\[
g:= W_{01}f + Kf +  \bigl(u f(-1) + v f(1) \bigr)\, \psi_1 + z\,
\psi_2 \in \cF_{\max}
\]
by \eqref{eqpW2} and Lemma~\ref{K}. Thus to prove that $w \in
\cF(A)$, it is enough to show that
\begin{equation} \label{eqfeq}
u\, g(-1)+ v\, g(1) = \eta_{11}\,[f,\psi_1]+\eta_{12}\,[f,\psi_2]+
\alpha\,\eta_{11}\,\bigl(u f(-1) +v  f(1) \bigr)
+\alpha\,\eta_{12}\, z.
\end{equation}
To this end we calculate
\begin{align*}
u\,g(-1) & = u\, \bigl( (W_{01}f)(-1) + (Kf)(-1) +
  \bigl(u f(-1) + v f(1) \bigr)\, \psi_1(-1) + z\, \psi_2(-1) \bigr) \\
 & = u\, (W_{01}f)(-1) + u\,(Kf)(-1)
  + \alpha\, |u|^2 \eta_{11}\, \bigl( u f(-1) + v f(1) \bigr)
    +  \alpha\, |u|^2\, \eta_{12}\, z
 \intertext{
from \eqref{eqpsib} and \eqref{eqdpsi}. Similarly
           }
v\,g(1) & = v\, \bigl( (W_{01}f)(1) + (Kf)(1)
   + \bigl( u f(-1) + v f(1) \bigr) \, \psi_1(1)
    + z\, \psi_2(1) \bigr) \\
 & =v\, (W_{01}f)(1)+ v\,(Kf)(1)
    +  \alpha \, |v|^2 \, \eta_{11} \, \bigl(u f(-1) + v f(1) \bigr)
     + \alpha \, |v|^2 \, \eta_{12} \, z .
\end{align*}
Adding and using  \eqref{eqpW}, \eqref{eqebcK} and \eqref{equv1}, we
obtain \eqref{eqfeq}. This completes the proof of the lemma.
\end{proof}

The theorem now follows from Theorem~\ref{W} and Lemma \ref{fin}.
\end{proof}

We now specialize Theorems~\ref{tk0sebc}, \ref{tn=k2} and
\ref{tn2k1} to some of our earlier examples. First we consider
Example~\ref{e35} (cf. Example~\ref{e42}).

\begin{corollary}
Assume that $p=1$ and $r$ is of order $\nu_0 > -1$ on a
half-neighborhood of $0$, and of order $\nu_1 > -1$ on both a right
half-neighborhood of $-1$ and a left half-neighborhood of $1$. Then
there is a basis for each root subspace of $A$, so that the union of
all these bases is a Riesz basis of $L_{2,|r|} \oplus
\nC_{|\Delta|}$.
\end{corollary}

Now we consider Examples~\ref{r42} and~\ref{e43}.

\begin{corollary} \label{cpero}
Assume that $p$ is even, $r$ is odd and that
Condition~{\rm~\ref{c0}} holds. If $k = 0$ or
Condition~{\rm~\ref{cat-1}} holds, then there is a basis for each
root subspace of $A$, so that the union of all these bases is a
Riesz basis of $L_{2,|r|} \oplus \nC_{|\Delta|}$.
\end{corollary}

As a simple illustration of this corollary we could consider the
eigenvalue problem stated in Example~\ref{esbe} but with $r$ odd and
of order $\nu_0$ at $0$ and  $\nu_1$ at $1$ (and hence of order
$\nu_1$ at $-1$, since $r$ is odd).

\begin{corollary}
Assume that $p$ is nearly even and $r$ is nearly odd.  If $k = 0$ or
Condition~{\rm~\ref{cat-1}} holds, then there is a basis for each
root subspace of $A$, so that the union of all these bases is a
Riesz basis of $L_{2,|r|} \oplus \nC_{|\Delta|}$.
\end{corollary}

\end{document}